\documentclass[12pt]{article}
\usepackage{amsfonts}
\usepackage{amssymb}
\usepackage{amsmath}
\usepackage{amsthm}
\usepackage{mathrsfs}
\def\trianglelefteqslant{\trianglelefteq}
\usepackage[latin1]{inputenc}
\usepackage[all]{xy}
\usepackage[left=1.20in,right=1.20in]{geometry}
\usepackage{xcolor}

\def\Z{\mathbb{Z}}

\def\dom{\backslash}
\def\normal{\mathop{\trianglelefteqslant}}

\newcommand{\ind}{{\rm Ind}}

\newcommand{\dfl}{{\rm Def }}
\newcommand{\ifl}{{\rm Inf }}
\newcommand{\indinf}{{\rm Indinf }}

\newcommand{\dtr}{{\rm det}}

\newcommand{\rac}{\mathbb Q}

\newcommand{\F}{\mathbb F}

\newcommand{\ent}{\mathbb Z}

\newcommand{\bizlie}[1]{\mathop{\,\raisebox{-.5ex}{$\widehat{\raisebox{.9ex}{\rule{2.5ex}{.07ex}}}_{#1}$}\,}}

\def\xar[#1]{\ar@{-}[#1]}
\def\sur{\overline}

\def\mpn{\medskip\par\noindent}
\def\smpn{\smallskip\par\noindent}
\def\mmpn{\vskip 1em minus 1em\par\noindent}

\def\sp{\bigskip\par}

\def\smp{\smallskip\par}

\def\pf{\par\smallskip\noindent{\bf Proof. }}

\def\endpf{~\hfill\rlap{\hspace{-1ex}\raisebox{.5ex}{\framebox[1ex]{}}\sp}\bigskip\pagebreak[3]}

\makeatletter
\renewenvironment{enumerate}{\ifnum \@enumdepth >3 \@toodeep\else
       \advance\@enumdepth \@ne
       \edef\@enumctr{enum\romannumeral\the\@enumdepth}\list
       {\csname  label\@enumctr\endcsname}{\setlength{\topsep}{1ex}
\setlength{\itemsep}{0 pt}\usecounter
         {\@enumctr}\def\makelabel##1{\hss\llap{##1}}}\fi}{\endlist}

\renewenvironment{itemize}{\ifnum \@itemdepth >3 \@toodeep\else
\advance\@itemdepth \@ne
\edef\@itemitem{labelitem\romannumeral\the\@itemdepth}
\list{\csname\@itemitem\endcsname}{\setlength{\topsep}{1ex}\setlength
{\itemsep}{0pt}\def\makelabel##1{\hss\llap{##1}}}\fi}
{\endlist}

\def\@seccntformat#1{\csname the#1\endcsname.\quad}

\def\section{\pagebreak[3]\setcounter{prop}{0}\setcounter{equation}{0}\@startsection{section}{1}{\z@}{4ex plus  6ex}{2ex}{\center\reset@font \large\bf}}
\newcommand{\subsect}[1]{\bigskip\par\noindent\pagebreak[3]\refstepcounter{subsection}\refstepcounter{prop}{\large\bf \thesection.\arabic{prop}.\ #1.\ }\bigskip\par}
\makeatother

\def\theprop{\thesection.\arabic{prop}}

\renewenvironment{equation}{\refstepcounter{subsection}\refstepcounter
{prop}$$}{\leqno{\bf (\theprop)}$$}

\newenvironment{rem}[1]{\refstepcounter{subsection}\refstepcounter
{prop} \mpn{{\bf \thesection.\arabic{prop}.}\ \ \bf#1.}}{\smp}
\newenvironment{enonce}[1]{\pagebreak[3]\refstepcounter{prop}\mmpn
{{\bf  \thesection.\arabic{prop}.\ #1.}}\begin{it} }{\end{it}\par}
\def\thesection{\arabic{section}}

\newenvironment{enonce*}[1]{\pagebreak[2]\smpn{#1}\begin{it} }{\end{it}\smp}

\newcommand{\result}[1]{\begin{enonce}{#1}}
\newcommand{\resultb}[1]{\begin{enonce*}{#1}}

\newcommand{\fresult}{\end{enonce}}
\newcommand{\fresultb}{\end{enonce*}}

\renewcommand{\leq}{\leqslant}
\renewcommand{\geq}{\geqslant}

\title{The Whitehead group of (almost) extra-special $p$-groups with $p$ odd}
\author{Serge Bouc\footnote{CNRS-LAMFA, UPJV, Amiens, France}~~and Nadia Romero\footnote{DEMAT, UGTO, Guanajuato, Mexico - Partially supported by SEP-PRODEP, project PTC-486}  
}

\def\rouge{}

\begin{document}
\centerline{\Large\bf The Whitehead group of }
\vspace{1ex}
\centerline{\Large\bf (almost) extra-special $p$-groups with $p$ odd}
\vspace{2ex}
\centerline{\bf Serge Bouc\footnote{CNRS-LAMFA, UPJV, Amiens, France}~{\rm and} Nadia Romero\footnote{DEMAT, UGTO, Guanajuato, Mexico - {Partially supported by SEP-PRODEP,  project PTC-486}}}
\vspace{4ex}
\begin{center}
\begin{minipage}{9.7cm}{
\begin{footnotesize}
{\bf Abstract:} Let $p$ be an odd prime number. We describe the Whitehead group of all extra-special and almost extra-special $p$-groups. For this we compute, for any finite $p$-group~$P$, the subgroup $Cl_1(\Z P)$ of $SK_1(\Z P)$, in terms of a genetic basis of $P$. We also introduce a deflation map $Cl_1(\Z P)\to Cl_1\big(\Z(P/N)\big)$, for a normal subgroup $N$ of $P$, and show that it is always surjective. {Along the way, we give a new proof of the result describing the structure of $SK_1(\Z P)$, when $P$ is an elementary abelian $p$-group.} \smallskip\par
{\bf Keywords:} Whitehead group, almost extra-special $p$-groups, genetic basis.\smallskip\par
{\bf MSC2010: } 19B28, 20C05, 20D15.
\end{footnotesize}
}
\end{minipage}
\end{center}
\section*{Introduction}
{Whitehead groups were introduced by J.H.C. Whitehead in~\cite{whitehead-simple}, as an algebraic continuation of {his work on} combinatorial homotopy.} The computation of the Whitehead group {$Wh(G)$} of a finite group {$G$} is in general very hard, and a compendium on the subject is the book by Bob Oliver (\cite{bob}) of 1988. Since then, {it seems that} not much progress has been made on this subject (see however~\cite{magu}, \cite{guaschi-juan-millan}, \cite{ushi}, \cite{ushitaki}).\par
Let $p$ be an odd prime number. In this paper, we describe the Whitehead group of all extra-special and almost extra-special $p$-groups. The main reason for focusing on these $p$-groups is that, apart from elementary abelian $p$-groups, they are exactly the finite $p$-groups all proper factor groups of which are elementary abelian (see e.g. Lemma~3.1 of~\cite{thecar}). In particular these groups appear naturally in various areas as first crucial step in inductive procedures (see \cite{thecar}, \cite{boma}, or the proof of Serre's Theorem in Section 4.7 of~\cite{benson2} for examples).\par
One of the main tools in our method is Theorem~9.5 of~\cite{bob}, which gives a first description, for a finite $p$-group $P$, of the subgroup $Cl_1(\Z P)$, an essential part of the torsion of $Wh(P)$. As this description requires the knowledge of the rational irreducible representations of~$P$, it seems natural to try to translate it in terms of a {\em genetic basis} of $P$, which provides an explicit description of the rational irreducible representations of~$P$, using only combinatorial informations on the poset of subgroups of $P$. The notions of genetic subgroup and genetic basis come from biset functor theory~(\cite{bouc}). They already have been used successfully for the computation of other groups related to representations of finite $p$-groups, e.g. the group of units of Burnside rings (see \cite{burnsideunits}), or the Dade group of endopermutation modules (see~\cite{dadegroup}). \par
Our first task is then to obtain a description of $Cl_1(\Z P)$ in terms of a genetic basis of $P$. {The advantage of t}his description {is that it} only requires data of a combinatorial nature, and makes no use of linear representations. In particular this makes it much easier to implement for computational purposes, using e.g. GAP software (\cite{GAP4}). {{This} also explains why in the proof of our main theorem (Theorem A below), we have to consider separately two special cases of ``small'' $p$-groups, for which the lattice of subgroups is in some sense too tight for the general argument to work.} \par
We observe next on this description that for any normal subgroup $N$ of $P$, there is an obvious {\em deflation operation} \mbox{$\dfl_{P/N}^P:Cl_1(\Z P)\to Cl_1\big(\Z(P/N)\big)$}. Even if there seems to be no inflation operation in the other direction, which would provide a section of this map, we show that $\dfl_{P/N}^P$ is always surjective. The existence of such a surjective deflation map already follows from Corollary 3.10 of~\cite{bob}, and it is shown in~\cite{K1-genotype} that our deflation map is indeed the same as the map defined there. \par
 As we will see in the next section, in the case of an extra-special or almost extra-special $p$-group~$P$, the group $Cl_1(\Z P)$ is equal to the torsion part $SK_1(\Z P)$ of $Wh(P)$, so the computation of the Whitehead group of $P$ comes down to the knowledge of $Cl_1\big(\Z\big(P/\Phi(P)\big)\big)$, where $\Phi(P)$ is the Frattini subgroup of $P$, and a detailed analysis of what happens with the unique faithful rational irreducible representation of $P$. Our main theorem is the following: 
\resultb{{\bf Theorem A.}} Let $p$ be an odd prime, and let $P$ be an extra-special $p$-group of order at least $p^5$ or an almost extra-special $p$-group of order at least $p^6$. Set $N=\Phi(P)=P'$.
\begin{enumerate}
\item The group $Cl_1(\Z P)$ is isomorphic to $K\times (C_p)^M$, where $M=\frac{p^{k-1}-1}{p-1}-\binom{p+k-2}{p}$ if $|P|=p^k$, and $K$ is the kernel of $\dfl_{P/N}^P:Cl_1(\Z P)\to Cl_1\big(\Z (P/N)\big)$.
\item {\rouge The group $K$ is cyclic. More precisely, it is:
\begin{enumerate}
\item trivial if $P$ is extra-special of order $p^5$ and exponent $p^2$.
\item of order $p$ if $P$ is almost extra-special of order $p^6$.
\item isomorphic to $Z(P)$ in all other cases.
\end{enumerate} }
\end{enumerate}

\fresultb

It should be noted that our methods not only give the algebraic structure of the Whitehead group, but allow also for the determination of explicit generators of its torsion subgroup {(see. Remark~\ref{generators} for details)}.\par
The paper is organized as follows: Section 1 is a review of definitions and basic results on Whitehead groups, genetic bases of $p$-groups, extra-special and almost extra-special $p$-groups. In particular, the subgroups $Cl_1(\Z G)$ and $SK_1(\Z G)$ of the Whitehead group $Wh(G)$ of a group $G$ are introduced. In Section 2, we give a procedure (Theorem~\ref{elmejor}) to compute $Cl_1(\Z P)$ for a finite $p$-group $P$ (with $p$ odd) in terms of a genetic basis of $P$. This procedure may be of independent interest, in particular from an algorithmic point of view. Finally, in Section~3, we begin by giving a new proof of the structure of $Wh(P)$ for an elementary abelian $p$-group $P$, for $p$ odd, and then we come to our main theorems. The first of them is  Theorem \ref{debut}, which will give us point 1 of Theorem A, then  Theorem \ref{principal} and Proposition \ref{fait chier} will deal with the description of the group $K$, giving us point 2. With these results we completely determine the structure of $Wh(P)$, when $P$ is an extra-special or almost extra-special $p$-group, for $p$ odd.
\section{Preliminaries}

Let $G$ be a group. We will write $Z(G)$ for the center of $G$ and $G'$ for its commutator subgroup. The Frattini subgroup of $G$ is denoted by $\Phi (G)$. If $H$ is a subgroup of $G$, the normalizer  of $H$ in $G$ will be denoted by $N_G(H)$, and its centralizer by $C_G(H)$.
If $H=\langle h\rangle$ for an element $h\in G$, we may also write $C_G(H)=C_G(h)$. \\

\subsect{About the Whitehead group}

Let $R$ be an associative ring with unit. The infinite general linear group of $R$, denoted by $GL(R)$, is defined as the direct limit of the inclusions $GL_n(R)\rightarrow GL_{n+1}(R)$ as the upper left block matrix. If, for each $n>0$, we denote by $E_n(R)$ the subgroup of $GL_n(R)$ generated by all elementary $n\times n$-matrices -- i.e. all those which are the identity except for one non-zero off-diagonal entry -- and we take the direct limit as before for the groups $E_n(R)$, we obtain a subgroup of $GL(R)$, denoted by $E(R)$. Whitehead's Lemma (Theorem 1.13 in \cite{bob}, for instance) states that $E(R)$ is equal to the derived subgroup of $GL(R)$. The group $K_1(R)$ is defined as the abelianization of $GL(R)$, which is then equal to $GL(R)/E(R)$.

If $G$ is a group and we take $R$ as the group ring $\ent G$, then elements of the form $\pm g$ for $g\in G$ can be regarded as invertible $1\times1$-matrices over $\ent G$ and hence they represent elements in $K_1(\ent G)$. Let $H$ be the subgroup of $K_1(\ent G)$ generated by classes of elements of the form $\pm g$ with $g\in G$. The {\em Whitehead group} of $G$ is defined as $Wh(G)=K_1(\ent G)/H$. 

If $G$ is finite, the groups $K_1(\ent G)$ and $Wh(G)$ are finitely generated abelian groups (see Theorem 2.5 in \cite{bob}).

For the rest of this section, $G$ denotes a finite group.

\result{Definition}
Let $D$ be an integral domain and $K$ be its field of fractions, then $SK_1(DG)$ denotes the kernel of the morphism
\begin{displaymath}
K_1(DG)\rightarrow K_1(KG).
\end{displaymath}
\fresult

By Theorem 7.4 in \cite{bob}, $SK_1(\ent G)$ is isomorphic to the torsion subgroup of $Wh(G)$. Hence, $Wh(G)$ is completely determined by $SK_1(\ent G)$ and the rank of its free part {(i.e. its {\em free rank})}. According to Theorem 2.6 in \cite{bob}, this free rank is equal to $r-q$, where $r$ is the number of non-isomorphic irreducible $\mathbb{R}$-representations of $G$ and $q$ is the number of non-isomorphic irreducible $\mathbb{Q}$-representations of $G$.

\result{Definition}
Consider the ring of $p$-adic integers $\hat{\ent}_p$. The group $Cl_1(\ent G)$ is defined as the kernel of the localization morphism
\begin{displaymath}
l:SK_1(\ent G)\rightarrow \bigoplus_p SK_1(\hat{\ent}_pG).
\end{displaymath}
\fresult
By Theorem 3.9 in \cite{bob}, $SK_1(\hat{\ent}_pG)$ is trivial whenever $p$ does not divide $|G|$, and $l$ is onto. In particular, $SK_1(\ent G)$ sits in an extension
\begin{displaymath}
0\longrightarrow Cl_1(\ent G)\longrightarrow SK_1(\ent G)\longrightarrow\bigoplus_p SK_1(\hat{\ent}_pG)\longrightarrow 0.
\end{displaymath}
This extension is used by Oliver in \cite{bob} to describe $SK_1(\ent G)$ in many cases. The important feature of the groups treated in this paper is that {their derived subgroups are central and so by Theorem 8.10 in  \cite{bob},} {the group} $\bigoplus_p SK_1(\hat{\ent}_pG)$ is trivial, {then} describing $SK_1(\ent G)$ amounts to describing $Cl_1(\ent G)$.

\vspace{2ex}
\pagebreak[3]
\subsect{About genetic bases}

A finite $p$-group $Q$ is called a {\em Roquette $p$-group} if it has normal $p$-rank 1, i.e. if all its abelian normal subgroups are cyclic. The Roquette $p$-groups (see \cite{roquette} or Theorem~4.10 of \cite{gorenstein} for details) of order $p^n$ are the cyclic groups $C_{p^n}$, if $p$ is odd, and the cyclic groups~$C_{2^n}$, the generalized quaternion groups $Q_{2^n}$, for $n\geq 3$, the dihedral groups~$D_{2^n}$, for $n\geq 4$, and the semidihedral groups $SD_{2^n}$, for $n\geq 4$, if $p=2$. A Roquette $p$-group~$Q$ admits a unique {faithful} rational irreducible representation $\Phi_Q$ (see e.g. Proposition~9.3.5 in~\cite{bouc}).

\result{Definition}
Let $P$ be a finite $p$-group. A subgroup $S$ of $P$ is called {\em genetic} if the section $S\trianglelefteqslant N_P(S)\leqslant P$ satisfies 
\begin{enumerate}
\item The group $N_P(S)/S$ is a Roquette group.
\item Let $\Phi=\Phi_{N_P(S)/S}$ be the only faithful irreducible $\rac$-representation of $N_P(S)/S$ and $V=\textnormal{Ind}_{\,N_P(S)}^{\,P}\textnormal{Inf}_{\,N_P(S)/S}^{\,N_P(S)}\phi$, then  the functor $\textnormal{Ind}_{\,N_P(S)}^{\,P}\textnormal{Inf}_{\,N_P(S)/S}^{\,N_P(S)}$ induces an isomorphism of $\rac$-algebras 
\[\textnormal{End}_{\rac P}V\cong \textnormal{End}_{\rac (N_P(S)/S)}\Phi.\]
\end{enumerate}
\fresult
Note that the right-hand side algebra is actually a skew field by Schur's lemma. So $\textnormal{End}_{\rac P}V$ is also a skew field, hence $V$ is an indecomposable -- that is, irreducible -- $\rac P$-module.

\result{Notation}
Let $P$ be a finite $p$-group and $S$ be a genetic subgroup of $P$.
We write
\begin{displaymath}
V(S)=\ind_{N_P(S)}^P\ifl_{N_P(S)/S}^{N_P(S)}\Phi_{N_P(S)/S}.
\end{displaymath}
\fresult
Then $V(S)$ is an irreducible $\rac$-representation of $P$. Conversely, by Roquette's Theorem (Theorem 9.4.1 in \cite{bouc}), for each irreducible $\rac$-representation $V$ of $P$, there exists a genetic subgroup $S$ of $P$ such that $V\cong V(S)$.

The following theorem characterizes the genetic subgroups of a $p$-group. First some notation: for a subgroup $S$ of a finite $p$-group $P$, let $Z_P(S)\geqslant S$ be the subgroup of $N_P(S)$ defined by $Z_P(S)/S=Z\big(N_P(S)/S\big)$. In particular $Z_P(S)=N_P(S)$ if $N_P(S)/S$ abelian, e.g. if $N_P(S)/S$ is a Roquette $p$-group for $p$ odd.

\result{Theorem}[Theorem 9.5.6 in \cite{bouc}]\label{cargen}
Let $P$ be a finite $p$-group and $S$ be a subgroup of $P$ such that $N_P(S)/S$ is a Roquette group. Then the following conditions are equivalent:
\begin{enumerate}
\item The subgroup $S$ is a genetic subgroup of $P$.
\item If $x\in P$ is such that $^xS\cap Z_P(S)\leqslant S$, then $^xS=S$.
\item If $x\in P$ is such that $^xS\cap Z_P(S)\leqslant S$ and $S^x\cap Z_P(S)\leqslant S$, then $^xS=S$.
\end{enumerate}
\fresult

The next result is part of Theorem 9.6.1 in \cite{bouc}.

\result{Theorem}\label{cariso}
Let $P$ be a finite $p$-group and $S$ and $T$ be genetic subgroups of $P$. The following conditions are equivalent:
\begin{enumerate}
\item The $\rac P$-modules $V(S)$ and $V(T)$ are isomorphic.
\item There exist $x,y\in P$ such that $^xT\cap Z_P(S)\leqslant S$ and $^yS\cap Z_P(T)\leqslant T$.
\item There exists $x\in P$ such that $^xT\cap Z_P(S)\leqslant S$ and $S^x\cap Z_P(T)\leqslant T$.
\end{enumerate}
If these conditions hold, then in particular the groups $N_P(S)/S$ and $N_P(T)/T$ are isomorphic.
\fresult

The relation between groups appearing in point 2 is denoted by $S\bizlie{P}T$. The theorem shows that this relation is an equivalence relation on the set of genetic subgroups of $P$, and we have the following definition.

\result{Definition}[Definition 9.6.11 in \cite{bouc}]
Let $P$ be a finite $p$-group. A {\em genetic basis} of~$P$ is a set of representatives of the equivalence classes of $\bizlie{P}$ in the set of genetic subgroups of $P$.
\fresult
\result{Lemma}\label{fidele}
Let $P$ be a finite $p$-group and $S$ be a genetic subgroup of $P$. 
\begin{enumerate}
\item The kernel of $V(S)$ is equal to the intersection of the conjugates of $S$ in $P$. 
\item In particular $V(S)$ is faithful if and only $S$ intersects $Z(P)$ trivially.
\end{enumerate}
\fresult
\pf
Denote by $\Phi$ the unique rational irreducible representation of the Roquette group $N_P(S)/S$. Then 
$$V(S)=\ind_{N_P(S)}^P\ifl_{N_P(S)/S}^{N_P(S)}\Phi\cong\bigoplus_{x\in[P/N_P(S)]}x\otimes \Phi,$$
where $[P/N_P(S)]$ is a chosen set of representatives of $N_P(S)$-cosets in $P$.
 An element $g\in P$ acts trivially on $V(S)$ if and only if it permutes trivially the summands of this decomposition, that is if $g^x\in N_{P}(S)$ for any $x\in {P}$, and if moreover $g^x$ acts trivially on $\Phi$, which means that $g^x\in S$, since $\Phi$ is faithful. This proves Assertion~1.\par
Assertion 2 follows, since
$$\big(\bigcap_{x\in P}{^xS}\big)\cap Z(P)=\bigcap_{x\in P}\big({^xS}\cap Z(P)\big)=\bigcap_{x\in P}{^x\big({S}\cap Z(P)\big)}=S\cap Z(P),$$
{and since the normal subgroup $T=\bigcap_{x\in P}\limits{^xS}$ of $P$ is trivial if and only if $T\cap Z(P)={\bf 1}$, i.e. if $S\cap Z(P)={\bf 1}$.}
\endpf

\begin{rem}{Remark}
If $P$ is abelian, then a subgroup $S$ of $P$ is genetic if and only if $P/S$ is cyclic. Moreover the relation $\bizlie{P}$ is the equality relation in this case, so there is a unique genetic basis of $P$, consisting of all the subgroups $S$ of $P$ such that $P/S$ is cyclic.
\end{rem}
\begin{rem}{Remark}\label{normal normalizer}
Let $P$ be a finite $p$-group, and $S$ be a subgroup of $P$ such that $N_P(S)$ is normal in $P$. If $N_P(S)/S$ is a Roquette group, then $S$ is a genetic subgroup of~$P$. Indeed if $N_P(S)$ is normal in $P$, then $N_P(S)=N_P({^xS})$ for any $x\in P$. Hence ${^xS}\normal N_P(S)$, and the group ${^x}S\cdot S/S$ is a normal subgroup of $N_P(S)/S$. It is trivial if and only if it intersects trivially the center $Z_P(S)/S$ of $N_P(S)/S$. Then
\begin{eqnarray*}
{^xS}=S\iff {^x}S\cdot S/S={\bf 1}&\iff& {^x}S\cdot S\cap Z_P(S)=\big({^xS}\cap Z_P(S)\big)S=S\\
&\iff& {^xS}\cap Z_P(S)\leqslant S.
\end{eqnarray*}
\end{rem}

\subsect{About extra-special and almost extra-special $p$-groups}

\result{Definition}
Let $p$ be a prime and $P$ be a finite $p$-group.
\begin{enumerate}
\item The group $P$ is called extra-special if $Z(P)=P'=\Phi(P)$ has order $p$.
\item The group $P$ is called almost extra-special if $P'=\Phi(P)$ has order $p$ and $Z(P)$ is cyclic of order $p^2$.
\end{enumerate}
\fresult
Extra-special and almost extra-special $p$-groups can be classified in the following way.

\result{Notation}
Let $H$, $K$ and $M$ be groups such that $M\leqslant Z(H)$ and such that there exists an injective map $\theta:M\rightarrow Z(K)$. The central product of $H$ and $K$ with respect to $\theta$ will be denoted by $H*_\theta K$, and simply by $H*K$ if $\theta$ is clear from the context.

For any integer $r\geq 1$, we will write $H^{*r}$ for the central product of $r$ copies of the group $H$, where $M=Z(H)$, with the convention $H^{*1}=H$.

For $p\neq 2$, set
\begin{displaymath}
M(p)=\langle x,\, y\mid x^p=y^p=1,\, [x^{-1},\, y]=[y,\,x]={}^y[y,\,x]\rangle
\end{displaymath}
and
\begin{displaymath}
N(p)=\langle x,\, y\mid x^{p^2}=y^p=1, {}^yx=x^{1+p}\rangle.
\end{displaymath}
\fresult

\result{Theorem}
\label{clasextra}
Let $p$ be a prime and $P$ be a finite $p$-group.
\begin{enumerate}
\item If $P$ is extra-special, then there exists an integer $r\geq 1$ such that $P$ has order $p^{2r+1}$ and $P$ is isomorphic to only one of the following groups: $D_8^{*r}$ or $Q_8*D_8^{*(r-1)}$ if $p=2$, and $M(p)^{*r}$ or $N(p)*M(p)^{*(r-1)}$ if $p\neq 2$.
\item If $P$ is almost extra-special, then there exists an integer $r\geq 1$ such that $P$ has order $p^{2r+2}$ and $P$ is isomorphic to only one of the following groups: $D_8^{*r}*C_4$ if $p=2$, and $M(p)^{*r}*C_{p^2}$ if $p\neq 2$.
\end{enumerate}
\fresult
\pf
The proof of 1 can be found in Section 5.5 of \cite{gorenstein}. As for point 2, one can refer to Sections 2 and 4 of \cite{thecar}.
\endpf

Observe that if $p$ is odd, the exponent of the group characterizes the isomorphism type of extra-special $p$-groups, one of them has exponent $p$ and the other one has exponent $p^2$.

If $P$ is an (almost) extra-special group, the quotient $P/P'$ is elementary abelian, so it can be regarded as a (finite-dimensional) vector space $E$ over the finite field $\F_p$. Moreover, if we take $z$ a generator of $P'$, then $E$ is endowed with a bilinear form 
\begin{displaymath}
b:E\times E\rightarrow \F_p,
\end{displaymath}
that sends an element $(u,\, v)$ to $b(u,\, v)$, the element of $\F_p$ satisfying $[\tilde{u},\, \tilde{v}]=z^{b(u,\, v)}$ for all $\tilde{u}\in u$, $\tilde{v}\in v$ and $u,\, v\in E$.  This bilinear form is  alternating, i.e. $b(v,\,v)=0$ for all $v\in E$, hence it is antisymmetric, i.e. $b(u,\,v)=-b(v,\,u)$ for all $u,v\in E$. 
Section 20 in \cite{soluble} is concerned with this bilinear form for extra-special groups, but the property of being alternating is called \textit{symplectic}. 

Recall that if $f:V\times V\rightarrow K$ is a bilinear form on a finite dimensional vector space~$V$ over a field $K$, its {\em left radical} $V^\perp$ is defined by $V^{\perp}=\{v\in V\mid f(v,\,w)=0\,\forall w\in V\}$, and its {\em right radical} by $^\perp V=\{v\in V\mid f(w,\,v)=0\,\forall w\in V\}$. Clearly $V^\perp={^\perp V}$ when $f$ is antisymmetric. The {\em rank} of $f$ is the codimension of $V^\perp$. The form $f$ is called {\em non-degenerate} if $V^\perp=\{0\}$.\par
With the help of Lemma 20.4 in \cite{soluble}, we have the following observation.

\result{Observation}
The bilinear form $b$ is non-degenerate if and only if $P$ is extra-special. If $P$ is almost extra-special, then $E^{\perp}=\pi(Z(P))$ is a line in $E$, where $\pi:P\rightarrow P/P'$ is the projection morphism.
\fresult

In section 3.2 we will use the following result, which is part of Lemma 2.6 in \cite{boma}. For the proof we refer the reader to this reference.

\result{Lemma}
\label{deboma}
Let $P$ be an (almost) extra-special $p$-group, and let $Q$ be a non-trivial subgroup of $P$. Then
\begin{enumerate}
\item $Q\trianglelefteqslant P\Leftrightarrow P'\leqslant Q$.
\item If $Q$ is not normal in $P$, then $N_P(Q)=C_P(Q)$. In particular, it follows that in this case $Q$ is elementary abelian of rank at most $r$, for the integer $r$ defined as in Theorem \ref{clasextra}, and we have $|Q||C_P(Q)|=|P|$. Moreover, $C_P(Q)=Q\times U$, where $U$ is (almost) extra-special of order $|P|/|Q|^2$ or $U=Z(P)$.
\end{enumerate}
\fresult

\section{$Cl_1$ of finite $p$-group algebras for odd $p$}

The goal of this section is to re-write Theorem 9.5 in \cite{bob} in terms of genetic bases, in the most possible succinct way. We take the statement of this theorem appearing in Section 1 of \cite{mine}, which says the following: let $p$ be an odd prime, let $P$ be a finite $p$-group, and write $\rac P\cong\prod_{i=1}^kA_i$, where $A_i$ is simple with irreducible module $V_i$ and center $K_i={\rm End}_{\rac P}V_i$. By Roquette's Theorem, for each $1\leq i\leq k$, the field $K_i$ is isomorphic to $\rac (\zeta_{r_i})$, where $\zeta_{r_i}$ is a primitive $p^{r_i}$-th root of unity for some non-negative~$r_i$, and $A_i$ is isomorphic to a matrix algebra over $\rac(\zeta_{r_i})$.

Consider the abelian group $T=\prod_{i=1}^k\limits\langle\zeta_{r_i}\rangle$. For each $h\in P$, define
\begin{displaymath}
\psi_h:C_P(h)\rightarrow T,\quad \psi_h(g)=(\dtr_{K_i}(g, V_i^h))_i
\end{displaymath}
where $V_i^h=\{x\in V_i\mid hx=x\}$. Here $V_i^h$ is viewed as a $K_iC_P(h)$-module, so $\dtr_{K_i}(g, V_i^h)$ is the determinant (in $K_i$) of the action of $g$ in $V_i^h$. Since $P$ is a $p$-group, this determinant is in $\langle \zeta_{r_i}\rangle$. Then Theorem 9.5 in \cite{bob} can be written in the following way.

\result{Theorem}\label{elbueno}
Let $p$ be an odd prime and consider $T$ and $\psi_h:C_P(h)\rightarrow T$, for each $h\in P$, as before. Then
\begin{displaymath}
Cl_1(\ent P)\cong T/\langle \textnormal{Im}\psi_h\mid h\in P\rangle.
\end{displaymath}
\fresult

Now, since $p$ is odd, if we let $\mathcal{S}=\{S_1,\ldots , S_k\}$ be a genetic basis of $P$, then $N_P(S_i)/S_i$ is cyclic for every $1\leq i\leq k$ and each simple $\rac P$-module $V_i$ is isomorphic to $V(S_i)=\indinf^P_{N_P(S_i)/S_i}\Phi_{N_P(S_i)/S_i}$. 
Then the abelian group $T$ defined before is isomorphic to $\Gamma(P)=\prod_{i=1}^k\limits \big(N_P(S_i)/S_i\big)$. This is because we can see the module $\Phi_{N_P(S_i)/S_i}\cong \rac (\zeta_{r_i})$, where $p^{r_i}=|N_P(S_i)/S_i|$ as actually being generated by a generator of $N_P(S_i)/S_i$ and thus the action of $N_P(S_i)/S_i$ on it can be seen as multiplication on the group. In particular, $\dtr_{K_i}(g, V_i^h)$, which is an element of the field $K_i$, can be regarded as an element in $N_P(S_i)/S_i$.  The first step in re-writing Theorem \ref{elbueno} is to find this element, for every $S_i$, every element $h\in P$ and every $g\in C_P(h)$. 
\result{Notation}
Let $p$ be and odd prime. If $V$ is a simple $\rac P$-module and \mbox{$S$} is a genetic subgroup of $P$ corresponding to $V$, we will write $\dtr_{N_P(S)/S}(g, V^h)$ for the element in $N_P(S)/S$ corresponding to 
$\dtr_{K}(g, V^h)$, where $K={\rm End}_{\rac P}(V)$. 
\fresult\mpn

\result{Lemma}\label{formula}
Suppose $p$ is an odd prime. Let $V$ be a simple $\rac P$-module and $S$ be a
genetic subgroup of $P$ corresponding to $V$. 
Take an element $h$ in~$P$ { and let $H=\langle h\rangle$. If $g$ is in $C_P (H)$, let}  $[H\langle g\rangle\backslash P/N_P(S)]$ be a set of representatives of the double cosets of $P$ on $H\langle g\rangle$ and $N_P(S)$.  
Then we have
{
\begin{displaymath}
\dtr_{N_P(S)/S}\big(g,\, V^h \big)=\prod_{\substack{x\in [H\langle g\rangle\backslash P/N_P(S)]\\s.t.\, H^x\cap N_P(S)\leqslant S}}\overline{l_{g,\, x}},
\end{displaymath}
where $\sur{l_{g,\, x}}$ is determined as follows: {the set} $I_{g,x}=\langle g\rangle\cap {\big(}H{\cdot}{^xN_P(S)}{\big)}$ is {actually} a subgroup of $\langle g\rangle$. Let $m=|\langle g\rangle : I_{g,x}|$. Then $g^m$ can be written {as} $g^m=h{\cdot}{^xl_{g,x}}$, for $h\in H$ and $l_{g,x}\in N_P(S)$; the  element $l_{g,\,x}$ may not be unique, but its class $\overline{l_{g,\, x}}$ in $N_P(S)/S$ is, thanks to the conditions on $x$.}
\fresult
\pf
Set $\Phi$ for $\Phi_{N_P(S)/S}$. Let $[P/N_P(S)]$ be a set of representatives of the cosets of~$P$ in $N_P(S)$. Since $V\cong \indinf^P_{N_P(S)/S}\Phi$, we can write it as $\mathop{\bigoplus}_{a\in[P/N_P (S)]}\limits a \otimes \Phi$. The action of $y \in P$ is given by $y(a \otimes \omega) = ya \otimes \omega$, which is equal to $\tau_y(a)\otimes \overline{n_{y,\, a}}\omega$, if $ya=\tau_y(a)n_{y,\,a}$ for a corresponding $n_{y,\, a}$ in $N_P (S)$, with $\overline{n_{y,\, a}}$ being its class in $N_P(S)/S$. 

If $y \in H$ fixes an element $u = \sum_{a\in[P/N_P (S)]}\limits a \otimes \omega_a$ of $V$, we have
\begin{displaymath}
\sum_{a\in [P/N_P(S)]}\tau_y(a) \otimes \overline{n_{y,\,a}}\omega_a =\sum_{a\in [P/N_P(S)]}a \otimes \omega_a=\sum_{a\in[P/N_P(S)]}\tau_y(a)\otimes \omega_{\tau_y(a)}. 
\end{displaymath}
That is, for every $a\in [P/N_P(S)]$ we should have $\overline{n_{y,\, a}}\omega_a = \omega_{\tau_y(a)}$. If $\tau_y(a)=a$, then we should have that 
$y$ is in $H\cap { }^a\!N_P (S)$ and  that $H^a\cap N_P (S) \leqslant S$, if $\omega_a$ is different from zero. We consider then the set $[H\backslash P/N_P(S)]$ and we have that
\begin{displaymath}
u=\sum_{\substack{x\in[H\backslash P/N_P(S)]\\ s.t.\, H^x\cap N_P(S)\leqslant S}}\sum_{z\in[H/H\cap { }^xN_P(S)] }zx\otimes \omega_x.
\end{displaymath}
This means that a $\rac$-basis for $V^h$ is given by $\mu_{x,\, \omega}=\sum_{z\in[H/H\cap { }^x N_P (S)]}\limits zx \otimes \omega$ with $x$ running over $F=\{x \in [H\backslash P/N_P (S)] \mid H^x \cap N_P (S)\leqslant S\}$ and $\omega$ running over a $\rac$-basis of $\Phi$. Since $\Phi\cong \rac(\zeta_r)\cong$End$_{\rac P}(V)$, where $\zeta_r$ is a primitive $p^r$-th root of unity, if the order of $N_P(S)/S$ is $p^r$, then a $\rac(\zeta_r)$-basis of $V^h$ is the set of $\mu_x=\mu_{x,\, 1}$ with $x$ running over $F$.

Now, for $\mu_x$ we have that $g\mu_x$ is equal to
\begin{displaymath}
\sum_{z\in[H/H\cap { }^x N_P (S)]} zh_{g,\, x}\sigma_g(x) \otimes \overline{n_{g,\, x}}\,1
\end{displaymath}
if $gx=h_{g,\,x}\sigma_g(x)n_{g,\,x}$. That is 
\begin{displaymath}
g\mu_{x,\, 1}=h_{g,\, x}\mu_{\sigma_g (x),\, \overline{n_{g,\, x}}}=\mu_{\sigma_g (x),\, \overline{n_{g,\, x}}},
\end{displaymath}
 since $\mu_{\sigma_g (x),\, \overline{n_{g,\, x}}}$ is in $V^h$. So we can write $g\mu_x=\overline{n_{g,\,x}}\mu_{\sigma_g(x)}$, with $\overline{n_{g,\,x}}$ seen as an element in the field $\rac(\zeta_r)$. This implies that the action of $g$ in $V^h$ is given by a
monomial matrix~$A$, the coefficient in the non-zero entry of a row being $\overline{n_{g, x}}$. Then, the
determinant of $A$ is the product of the signature of the permutation $\sigma_g$ by the product of the coefficients $\overline{n_{g,\, x}}$. Since $p$ is odd, the signature is $+1$, and 
$$ \det(A)=\prod_{\substack{x\in [H\backslash P/N_P(S)]\\s.t.\, H^x\cap N_P(S)\leqslant S}}\overline{n_{g,\, x}}.
 $$
{
Observe now that the group $H^x\cap N_P(S)\leq S$ does not change if we replace $x$ by $yx$, for $y\in\langle g\rangle$, since $g$ centralizes $H$. Moreover the intersection $I_{g,x}=\langle g\rangle\cap {\big(}H{\cdot}{^xN_P(S)}{\big)}$ is equal to $\langle g\rangle\cap {\Big(}H{\cdot}\big({^xN_P(S)}\cap C_P(H\langle g\rangle)\big){\Big)}$, by Dedekind's modular law {\rouge(\cite{huppert}, Hilfsatz 2.12.c),} since $g$ and $H$ centralize $H$. Hence $I_{g,x}$ is a subgroup of $\langle g\rangle$. Now a set $R$ of representatives of $H\dom P/N_P(S)$ is the set of elements $yx$, where $x$ is in a set of representatives of $H\langle g\rangle\backslash P/N_P(S)$, and $y$ is in a set of representatives of $\langle g\rangle/I_{g,x}$. The set of elements $y$ can be taken as $\{1,g,g^2,\ldots,g^{m-1}\}$. For $y=g^i$, with $0\leq i\leq m-2$, we have $gyx=g^{i+1}x\in R$, so $n_{g,yx}=1$. {F}or $i=m-1$, we have $gyx=g^mx=h\cdot x\cdot l_{g,x}$.\endpf
}

Our final version of Theorem \ref{elbueno} is the following.

\def\Q{\mathbb{Q}}
\result{Theorem}\label{elmejor} Let $p$ be an odd prime and $P$ be a finite $p$-group. Take a set $\mathcal{C}$ of representatives of conjugacy classes of cyclic subgroups of $P$. For each $H\in\mathcal{C}$, let $\sur{E}_H$ be a generating set of the factor group $C_P(H)/H$ and $E_H\subseteq C_P(H)$ be a set of representatives of the classes $gH\in \sur{E}_H$. Let also $\mathcal{S}$ be a genetic basis of $P$ and for each $S\in \mathcal{S}$, let $[H\langle g\rangle\backslash P/N_P(S)]$ be a set of representatives of the double cosets of $P$ on $H\langle g\rangle$ and $N_P(S)$. Then 
$$Cl_1(\ent P)\cong \left(\prod_{S\in\mathcal{S}}\big(N_P(S)/S\big)\right)\Big/\mathcal{R},$$
where $\mathcal{R}$ is the subgroup generated by the elements $u_{H,g}=(u_{H,g,S})_{S\in\mathcal{S}}$, for $H\in\mathcal{C}$ and $g\in E_H$; {for $S\in\mathcal{S}$, the component $u_{H,g,S}$ of $u_{H,g}$ is given by
\begin{equation}\label{formule2}
u_{H,g,S}=\prod_{\substack{x\in [H\langle g\rangle\backslash P/N_P(S)]\\\rule{0ex}{1.3ex}s.t.\, H^x\cap N_P(S)\leqslant S}}\overline{l_{g,\, x}},
\end{equation}

\noindent where $\sur{l_{g,x}}$ is the image in $N_P(S)/S$ of the element $l_{g,x}\in N_P(S)$ determined as follows: let $m$ denote the index {of $\langle g\rangle\cap \big(H{\cdot}{^xN_P(S)}\big)$ in $\langle g\rangle$}. Then $g^m\in H{\cdot}{^xN_P(S)}$, so {$g^m$ can be written as} $g^m=h.{^xl_{g,x}}$, for some $h\in H$ and $l_{g,x}\in N_P(S)$.
}
\fresult
\pf As we said at the beginning of the section, since $p$ is odd, we have $Cl_1(\ent P)\cong\Gamma(P)/\mathcal{R}$, where
$$\Gamma(P)=\prod_{S\in\mathcal{S}}\big(N_P(S)/S\big)$$
and $\mathcal{R}$ is the subgroup generated by all the elements $u_{h,g}=(u_{h,g,S})_{S\in\mathcal{S}}$,  with $g\in C_P(h)$ and $u_{h,g,S}=\det_{N_P(S)/S}\big(g,V(S)^h\big)$, where
$V(S)= \indinf_{N_P(S)/S}^P\Phi_{N_P(S)/S}$. 

We first observe that $u_{h,g}=u_{^yh,^yg}$, for any $y\in G$. Indeed, setting $V=V(S)$, we have a commutative diagram
$$\xymatrix{
V^h\ar[r]^-g\ar[d]_-y&V^h\ar[d]^-y\\
V^{^yh}\ar[r]^-{^yg}&V^{^yh}
}
$$
where the arrows are given by the actions of the labelling elements. It follows that the determinant of $^yg$ acting on $V^{^yh}$ is equal to the determinant of $g$ acting on $V^h$. Hence to generate the subgroup $\mathcal{R}$ of $\Gamma(P)$, it suffices to take the elements $u_{h,g}$, where $(h,g)$ runs through a set of representatives of conjugacy classes of pairs of commuting elements in $P$.\par
Now clearly for each $h\in H$, the map $g\mapsto u_{h,g}$ is a group homomorphism from $C_P(h)$ to $\Gamma(P)$, hence $\mathcal{R}$ is generated by the elements $u_{h,g}$, where $h\in P$ and $g$ runs through a set of generators of $C_P(h)$. Moreover, as $h$ acts as the identity on $V^h$, the group generated by $h$ is contained in the kernel of this morphism.\par
Finally, by Lemma~\ref{formula}, setting $H=\langle h\rangle$, we have
$$u_{h,g,S}=\prod_{\substack{x\in [H\langle g\rangle\backslash P/N_P(S)]\\s.t.\, H^x\cap N_P(S)\leqslant S}}\overline{l_{g,\, x}}$$
and this depends only on $H$, so we may denote it by $u_{H,g,S}$, and by $u_{H,g}$ the corresponding element of $\Gamma(P)$. \par
It follows that $\mathcal{R}$ is generated by the elements $u_{H,g}$, where $H$ is a cyclic subgroup of~$P$ up to conjugation, and for a given $H$, the element $g$ runs through a subset of $C_P(H)$ which, together with $H$, generates $C_P(H)$. {Together with Lemma~\ref{elbueno}, this completes the proof.}
\endpf

\vspace{-2ex}
To finish the section, we {observe} that if $N$ is a normal subgroup of a finite \mbox{$p$-group~$P$} with $p$ odd, then there is surjective deflation morphism 
$$\dfl_{P/N}^P:Cl_1(\ent P)\to Cl_1\big(\ent (P/N)\big).\vspace{-2ex}$$
\result{Proposition}\label{ladef} Let $p$ be an odd prime and $P$ be a finite $p$-group. Suppose $N$ is a normal subgroup of $P$. Let $\mathcal{S}$ be a genetic basis of $P$ and $\mathcal{S}_N=\{S\in\mathcal{S}\mid N\leqslant S\}$. Let $\tilde{B}$ be the set of subgroups $\tilde{S}=S/N$ of $\tilde{P}=P/N$, for $S\in\mathcal{S}_N$. \begin{enumerate}
\item The set $\{\tilde{S}\mid S\in\mathcal{S}_N\}$ is a genetic basis of $\tilde{P}$, and for $S\in\mathcal{S}_N$, the projection $P\mapsto \tilde{P}$ induces an isomorphism $\pi_S:N_P(S)/S\to N_{\tilde{P}}(\tilde{S})/\tilde{S}$.
\item The composition
\[\xymatrix{
*!U(0.3){s:\Gamma(P)=\prod_{S\in\mathcal{S}}\limits\big(N_P(S)/S\big)}\ar[r]&*!U(0.3){ \prod_{S\in\mathcal{S}_N}\limits\big(N_P(S)/S\big)}\ar[r]^-{\prod\pi_S}&*!U(0.3){\prod_{\tilde{S}\in\tilde{\mathcal{S}}}\limits\big(N_{\tilde{P}}(\tilde{S})/\tilde{S}\big)=\Gamma({\tilde{P}})}
}
\]
induces a surjective deflation morphism $\dfl_{P/N}^P:Cl_1(\ent P)\to Cl_1\big(\ent (P/N)\big)$. {In particular $Cl_1\big(\Z (P/N)\big)$ is isomorphic to a quotient of $Cl_1(\Z P)$.}
\end{enumerate}
\fresult
\pf For Assertion 1, it is clear from the definitions that if $S/N$ is a genetic subgroup of $P/N$, then $S$ is a genetic subgroup of $P$. Moreover the irreducible representation of $P$ associated to $S$ is obtained by inflation from $P/N$ to $P$ of the irreducible representation of $P/N$ associated to $S/N$, up to the obvious isomorphism $\pi_S: N_P(S)/S\cong N_{\tilde{P}}(\tilde{S})/\tilde{S}$.\par
For Assertion 2, as the map $s:\Gamma(P)\to\Gamma({\tilde P})$ is surjective, all we have to check is that the subgroup $\mathcal{R}$ of defining relations for $Cl_1(\ent P)\cong\Gamma(P)/\mathcal{R}$ is mapped by $s$ inside the corresponding subgroup $\tilde{\mathcal{R}}$ of defining relations for $Cl_1(\ent \tilde{P})\cong\Gamma({\tilde{P}})/\tilde{\mathcal{R}}$. So let $H$ be a cyclic subgroup of $P$, let $g\in C_P(H)$, and let $S\in\mathcal{S}_N$. Then $\tilde{H}=HN/N$ is cyclic, and $\tilde{g}=gN\in C_{\tilde{P}}(\tilde{H})$. Moreover the map
$$H\dom P/N_P(S)\ni HxN_P(S)\mapsto \tilde{H}\tilde{x}N_{\tilde{P}}(\tilde{S})\in \tilde{H}\dom\tilde{P}/N_{\tilde{P}}(\tilde{S}),$$
where $\tilde{x}=xN\in\tilde{P}$, is a bijection, since $HxN_P(S)=HxNN_P(S)=HNxN_P(S)$ as $N_P(S)\geq S\geq N$. Hence we may identify the sets of representatives $[H\dom P/N_P(S)]$ and $[\tilde{H}\dom\tilde{P}/N_{\tilde{P}}(\tilde{S})]$ via this map. Moreover
$$H^x\cap N_P(S)\leq S\iff (HN)^x\cap N_P(S)\leq S\iff\tilde{H}^{\tilde{x}}\cap N_{\tilde{P}}(\tilde{S})\leq \tilde{S}.$$
Now, for $x\in  [H\dom P/N_P(S)]$ such that $H^x\cap N_P(S)\leq S$, the equality $gx=h_{g,x}\sigma_g(x)n_{g,x}$, where $h_{g,x}\in H$, $\sigma_g(x)\in [H\dom P/N_P(S)]$, and $n_{g,x}\in N_P(S)$, implies $\tilde{g}\tilde{x}=\widetilde{h_{g,x}}\widetilde{\sigma_g(x)}\widetilde{n_{g,x}}$ in $P/N$. In other words $\pi_S\big(\sur{n_{g,x}}\big)=\sur{n_{\tilde{g},\tilde{x}}}$, that is $s(u_{H,g})=u_{\tilde{H},\tilde{g}}$. This completes the proof.
\endpf

\begin{rem}{Remark} It is shown in \cite{K1-genotype} (where $\Gamma(P)$ is called the {\em genome} of $P$) that the map $s$ is the deflation map in a structure of {\em biset functor} on $\Gamma$, but we will not need the corresponding additional operations of induction, restriction and inflation in this paper. {It also follows from \cite{K1-genotype} that the deflation map $\dfl_{P/N}^P:Cl_1(\Z P)\to Cl_1\big(\Z(P/N)\big)$ we obtain here is the same as the map given by Corollary~3.10 of~\cite{bob}.}
\end{rem}

\section{Computing some Whitehead groups}

As we noted in Section 1, the examples we will consider in this section are all finite $p$-groups with $p$ odd that satisfy that the group $SK_1(\ent P)$ is equal to $Cl_1(\Z P)$. Hence, we will use Theorem \ref{elmejor} to calculate $Cl_1(\Z P)$. If $P$ is abelian, Theorem~\ref{elmejor} has a simpler expression, as it was already noted in Observation 1.13 of~\cite{mine}.

As we said in the introduction, to calculate the free rank of the Whitehead group of the groups in question, we will use Theorem 2.6 in \cite{bob}. We will also use Exercise~13.9 in Serre \cite{serre}, which says that if $G$ is a group of odd order and $c$ is the number of irreducible non-isomorphic complex representations of $G$, then $(c+1)/2$ is the number of irreducible non-isomorphic real representations of $G$.

We introduce some notation that will be helpful in both of our examples.

\result{Notation}\label{symmetric}
Let $p$ be a prime. Suppose that $W$ is a finite-dimensional vector space over the finite field $\F_p$.
We denote by $S(W)$ the symmetric algebra of $W$ and by $S^p(W)$ its homogeneous part of degree $p$.  If $\psi:W\rightarrow \F_p$ is a linear functional, the map 
$$w_1\otimes_{\F_p}\ldots\otimes_{\F_p}w_n \in W^{\otimes n} \mapsto \psi(w_1)\ldots\psi(w_n)\in\F_p$$
induces a well defined linear functional $S(W)\to\F_p$, that we denote by $A\mapsto A(\psi)$.\par
The choice of a basis $\{x_1,\ldots,x_k\}$ of $W$ over $\F_p$ yields a standard identification of $S(W)$ with the polynomial ring $\F_p[x_1,\ldots,x_k]$.\fresult
With such an identification, if $A=A(x_1,\ldots,x_k)\in S(W)$ and $\psi$ is a linear form on~$W$, then $A(\psi)=A\big(\psi(x_1),\ldots,\psi(x_k)\big)$. In particular $A(\psi)=0$ for all $\psi$ if and only if the polynomial function associated to $A$ is equal to zero, that is if $A(r_1,\ldots,r_k)=0$ for any $(r_1,\ldots,r_k)\in\F_p^k$: indeed since  $\{x_1,\ldots,x_k\}$ is a basis of $W$, for any such $k$-tuple $(r_1,\ldots,r_k)\in\F_p^k$, there is a unique linear form $\psi$ on~$W$ such that $\psi(x_i)=r_i$ for $1\leq i\leq k$.

\subsect{Elementary abelian $p$-groups}

\result{Lemma}
Let $p$ be an odd prime and $P$ be an elementary abelian $p$-group of rank~$k$, the free rank of $Wh(P)$ is equal to
{
\begin{displaymath}
\frac{(p^k-1)(p-3)}{2(p-1)}.
\end{displaymath}
}
\fresult
\pf
The number of non-isomorphic irreducible $\mathbb{R}$-representations of $P$ is $(p^k+1)/2$. On the other hand, the number of non-isomorphic irreducible $\mathbb{Q}$-representations of $P$ is equal to $(p^k+p-2)/(p-1)$, since the genetic basis for $P$ is given by all its subgroups of index $p$ plus $P$ itself. The result follows from Theorem 2.6 in \cite{bob}.
\endpf

The description of $SK_1$ for elementary abelian groups appeared first in Alperin et al. \cite{ados}. {We prove this result using our combinatorial approach with genetic bases. Our proof has some similarities with the one in [1], but makes no use of characters. It will also be useful when dealing with extra-special $p$-groups.} 

\result{Lemma} \label{Bxy generate}
\label{primero}
Let $p$ be a prime and let $W$ be a finite-dimensional vector space over~$\F_p$.
For $x$ and $y$ in $W$, we set $B_{x,\, y}=x^{p-1}y\in S^p(W)$. Then $S^p(W)$ is generated by the elements $B_{x,\, y}$ with $x$ and $y$ running over $W$.
\fresult
\pf
Recall that if $x_1,\ldots ,x_m$ are $m$ (not necessarily different) commuting variables, then for any $n$
\begin{displaymath}
(x_1+\cdots +x_m)^n=\sum_{\substack{\alpha_1,\ldots ,\alpha_m\textnormal{ s.t. }\\ \alpha_1+\cdots +\alpha_m=n}} \frac{n!}{\alpha_1!\cdots \alpha_m!}x_1^{\alpha_1}\cdots x_m^{\alpha_m}.
\end{displaymath}
This allows us to show that for any $n$
\begin{equation}\label{the formula}
\sum_{\emptyset \neq A\subseteq \{1,\ldots ,n\}}(-1)^{n-|A|}\left(\sum_{i\in A}x_i\right)^n=n!\, x_1\cdots x_n,
\end{equation}
and so if $p$ is prime, then
\begin{displaymath}
\sum_{\emptyset \neq A\subseteq \{1,\ldots ,p-1\}}(-1)^{p-1-|A|}\left(\sum_{i\in A}x_i\right)^{p-1}x_p=(p-1)!\, x_1\cdots x_p,
\end{displaymath}
which by Wilson's Lemma gives
\begin{displaymath}\sum_{\emptyset \neq A\subseteq \{1,\ldots ,p-1\}}(-1)^{p-|A|}\left(\sum_{i\in A}x_i\right)^{p-1}x_p=x_1\cdots x_p,
\end{displaymath}
that is
\begin{equation}\label{the formula 2}
\sum_{\emptyset \neq A\subseteq \{1,\ldots ,p-1\}}(-1)^{p-|A|}B_{\sum_{i\in A}\limits x_i,x_p}=x_1\cdots x_p.
\end{equation}
This completes the proof, by taking $x_1,\ldots,x_p\in W=S^1(W)$, which indeed commute in $S(W)$.
\endpf

\result{Theorem}[Theorem 2.4 in \cite{ados}]\label{abelem}
Let $p$ be an odd prime and $P$ be an elementary abelian $p$-group of rank $k$, then $SK_1(\ent P)$ is isomorphic to $(C_p)^N$ where
\begin{displaymath}
N=\frac{p^k-1}{p-1}-\binom{p+k-1}{p}.
\end{displaymath}
\fresult
\pf
As we said before, the genetic basis of $P$ consists of $P$ itself and all its subgroups of index $p$, so $SK_1(\ent P)$ is isomorphic to the quotient of
\begin{displaymath}
\Gamma(P)=\prod_{[P:Q]=p}(P/Q)
\end{displaymath}
by the subgroup generated by the elements $u_{x,\, y}$ for $x$, $y$ in $P$, where
\begin{displaymath}
(u_{x,\, y})_Q=\left\{\begin{array}{cl}
yQ & \textrm{if }x\in Q\\
1 & \textrm{otherwise.}
\end{array}\right.
\end{displaymath}
On the other hand, we can see $P$ as a vector space over $\mathbb{F}_p$, and for each subgroup~$Q$ of index $p$, i.e. for each hyperplane $Q$ of $P$, consider $\psi_Q:P\rightarrow \F_p$, a linear functional with kernel $Q$. The product of these $\psi_Q$ induces an isomorphism from $\Gamma(P)$ to 
\begin{displaymath}
V=\prod_{[P:Q]=p}\mathbb{F}_p,
\end{displaymath}
and the elements $u_{x,\, y}$ can be seen as
\begin{displaymath}
(u_{x,\, y})_Q=\left\{\begin{array}{cl}
\psi_Q(y) & \textrm{if }x\in Q\\
0 & \textrm{otherwise.}
\end{array}\right.
\end{displaymath}
We define a morphism $r:S^p(P)\rightarrow V$, sending $A\in S^p(P)$ to the vector whose $Q$-component is equal to $A(\psi_Q)$. 
We will show that Im$(r)$ is equal to the subspace of~$V$ generated by the elements $u_{x,\, y}$, and that $r$ is injective. This will give us the result. 

We show first that the elements $u_{x,\, y}$ are in Im$(r)$. For $x,\,y \in P$, let $B_{x,\, y}=x^{p-1}y\in S^p(P)$. Then
\begin{displaymath}
B_{x,\, y}(\psi_Q)=\psi_Q(x)^{p-1}\psi_Q(y)=\left\{\begin{array}{cl}
0 & \textrm{if }x\in Q\\
\psi_Q(y) & \textrm{otherwise}
\end{array}\right.
\end{displaymath}
since $\lambda^{p-1}=1$ if $\lambda$ is in $\mathbb{F}_p-\{0\}$ and $0^{p-1}=0$. In particular $B_{y,y}(\psi_Q)=\psi_Q(y)$ for all $y\in P$, thus $r(B_{y,y}-B_{x,\, y})=u_{x,\,y}$. \par
On the other hand, by Lemma~\ref{Bxy generate}, $S^p(P)$ is generated by the elements $B_{x,\, y}$ where $x$ and $y$ run through $P$, so we have that Im$(r)$ is contained in (hence equal to) the subspace of $V$ generated by the elements $u_{x,\,y}$. 

Finally, we prove that $r$ is injective. Let $A$ be in the kernel of $r$, then $A(\psi_Q)=0$ for every $Q$ of index $p$ in $P$. If $\psi$ is any other linear functional of $P$ with kernel $Q$, then there exists $\lambda \in \mathbb{F}_p$ such that $\psi = \lambda\psi_Q$, and so $A(\psi)=\lambda^pA(\psi_Q)=0$, since $A$ is homogeneous of degree $p$. Choosing a basis of $P$ over $\F_p$ as in Notation~\ref{symmetric}, we can view $A$ as a homogeneous polynomial of degree $p$, and the  polynomial function $A$ is zero. It remains to see that $A$ is actually the zero polynomial, but this {follows from Lemma~2.1 of~\cite{ados}.}
\endpf

{
\result{Corollary} Let $p$ be an odd prime, and $P$ be a finite $p$-group. If $|P/\Phi(P)|=p^k$, then $SK_1(\Z P)$ has a subquotient isomorphic to $(C_p)^N$, where 
\begin{displaymath}
N=\frac{p^k-1}{p-1}-\binom{p+k-1}{p},
\end{displaymath}
and in particular $SK_1(\Z P)\neq 0$ if $k\geq 3$.
\fresult
}
\pf Indeed $Cl_1\big(\Z\big(P/\Phi(P)\big)\big)\cong (C_p)^N$, for $N=\frac{p^k-1}{p-1}-\binom{p+k-1}{p}$, by Theorem~\ref{abelem}. Moreover, this group is a quotient of $Cl_1(\Z P)$, by Proposition~\ref{ladef}. Finally $Cl_1(\Z P)$ is a subgroup of $SK_1(\Z P)$, and $N>0$ if $k\geq 3$.\endpf

\subsect{Extra-special and almost extra-special $p$-groups}

We begin by finding a genetic basis of an (almost) extra-special $p$-group. 
\result{Proposition}\label{genetic basis extra-special}
Let $p$ be a prime and $P$ be an (almost) extra-special $p$-group. A genetic basis of $P$ is given by all its subgroups of index $p$, together with $P$ and a subgroup~$Y$ of maximal order such that $Y\cap Z(P)={\bf 1}$. In particular $P$ has a unique faithful rational irreducible representation, up to isomorphism.
\fresult
\pf
We abbreviate $Z(P)$ by $Z$. 

By theorems \ref{cargen} and \ref{cariso}, the subgroups of $P$ of index 1 or $p$ are genetic and are not linked modulo $\bizlie{P}$. They clearly intersect $Z$ non-trivially. On the other hand, any genetic subgroup $S\neq P$ which intersects $Z$ non-trivially must have index $p$, since $P'\leqslant S$ and so the cyclic group $P/S\cong (P/P')/(S/P')$, should have order $p$. 
This implies that if there is another group in $\mathcal{S}$, it must intersect $Z$ trivially. 

Let $Y\leqslant P$ be of maximal order with the property $Y\cap Z={\bf 1}$. By Lemma \ref{deboma}, we have that $C_P(Y)=N_P(Y)=YZ$. In particular $N_P(Y)\normal P$ and $N_P(Y)/Y\cong Z$ is a Roquette group. Then $Y$ is genetic, by Remark~\ref{normal normalizer}. Also, by  Lemma \ref{deboma}, we have that if $Y_1\leqslant P$ is a group such that $Y_1\cap Z={\bf 1}$, but it is not maximal order with this property, then $C_P(Y_1)=N_P(Y_1)$, but $N_P(Y_1)/Y_1$ is not cyclic. Thus $Y_1$ is not a genetic subgroup of $P$.

Finally, if $Y$ is a subgroup of $P$ of maximal order such that $Y\cap Z={\bf 1}$, then $Y$ has $|P/YZ|=|Y|$ distinct conjugates in $P$, by Lemma~\ref{deboma}. These conjugates are subgroups of index $p$ in the elementary abelian group $YP'$, and they all intersect trivially (that is, they don't contain) the group $P'$. Since there are exactly $|Y|$ subgroups of $YP'$ not containing $P'$, these subgroups are exactly the conjugates of $Y$ in $P$.\par
Now if $Y_0$ is another subgroup of $P$ such that $Y_0\cap Z={\bf 1}$, then $Y_0\cap YP'$ is a subgroup of $YP'$ which does not contain $P'$. Hence it is contained in some conjugate of~$Y$, and there exists $x\in P$ such that $Y_0\cap YP'\leqslant Y^x$. It follows that $Y_0\cap YZ\leqslant Y^x$, for $YP'$ is the subgroup of $YZ$ consisting of elements of order at most $p$. In other words 
$${^xY_0}\cap Z_P(Y)={^xY_0}\cap YZ={^x(Y_0\cap YZ)}\leqslant Y.$$
Now if $Y_0$ is another subgroup of maximal order such that $Y_0\cap Z={\bf 1}$, exchanging the roles of $Y$ and $Y_0$ in the previous argument shows that there also exists an element $y\in P$ such that ${^yY}\cap Y_0Z\leqslant Y_0$. By Theorem~\ref{cariso}, it follows that $Y_0\bizlie{P}Y$. The last assertion now follows from Lemma~\ref{fidele}.
\endpf

As a first consequence of this result we have.

\result{Lemma} Let $p$ be an odd prime, and $n$ be a positive integer. 
\begin{enumerate}
\item Let $P$ be an extra-special $p$-group of order $p^{2n+1}$. Then the free rank of $Wh(P)$ is equal to
{
$$\frac{(p^{2n}+p-2)(p-3)}{2(p-1)}.$$
}
\item Let $P$ be an almost extra-special group of order $p^{2n+2}$. Then the free rank of $Wh(P)$ is equal to
{
$$\frac{(p^{2n+1}+p^2+p+1)(p-3)+8}{2(p-1)}.$$
}
\end{enumerate}
\fresult
\pf The free rank of $Wh(P)$ is equal to $r-q$, where $r$ (resp. $q$) is the number of irreducible real (resp. rational) representations of $P$, up to isomorphism. In general, for a finite $p$-group $P$ and a field $K$ of characteristic 0, the irreducible representations of $P$ over $K$ can be recovered from the knowledge of a genetic basis $\mathcal{B}$ of $P$ {(see~\cite{barker-genotype}):  this is because the functor $R_K$ of representations of $p$-groups over $K$ is {\em rational} in the sense of Definition~10.1.3 of~\cite{bouc}, as can be easily deduced from Theorem~10.6.1 of~\cite{bouc}. A proof of this fact can also be found in~\cite{fast-roquette}.} In particular, the number $l_K(P)$ of such representations, up to isomorphism, is equal to
$$l_K(P)=\sum_{S\in\mathcal{B}}\partial l_K\big(N_P(S)/S\big),$$
where $\partial l_K(Q)$ denotes the number of {\em faithful} irreducible representations of a $p$-group $Q$ over $K$, up to isomorphism. For a Roquette $p$-group $Q$, we have moreover $\partial l_K(Q)=1$ if $Q={\bf 1}$, and $\partial l_K(Q)=l_K(Q)-l_K(Q/Z)$ otherwise, where $Z$ is the unique central subgroup of order $p$ of $Q$.\par
If $p$ is odd, all the groups $N_P(S)/S$, for $S\in\mathcal{B}$, are cyclic. Now for $Q=C_{p^m}$, with $m\geq 0$, we have 
$$l_{\mathbb{R}}(Q)=\frac{p^m+1}{2}\;\;\hbox{and}\;\;l_{\rac}(Q)=m+1.$$
It follows that $\partial l_{\mathbb{R}}(Q)=\frac{p^m-p^{m-1}}{2}$ if $m>0$, and $\partial l_{\mathbb{R}}(Q)=1$ if $m=0$. On the other hand $\partial l_{\rac}(Q)=1$ for any $m$. \par
In case $P$ is extra-special of order $p^{2n+1}$, the genetic basis obtained in Proposition~\ref{genetic basis extra-special} consists of the group $S=P$, for which $N_P(S)/S$ is trivial, of $\frac{p^{2n}-1}{p-1}$ subgroups $S$ of index $p$ in $P$, for which $N_P(S)/S\cong C_p$, and the subgroup $S=Y$, for which $N_P(S)/S\cong Z(P)\cong C_p$. This gives
$$r=l_{\mathbb{R}}(P)=1+ \frac{p^{2n}-1}{p-1}\frac{p-1}{2}+\frac{p-1}{2}=\frac{p^{2n}+p}{2},$$ 
and
$$q=l_{\rac}(P)=1+\frac{p^{2n}-1}{p-1}+1=\frac{p^{2n}+2p-3}{p-1}.$$
In case $P$ is almost extra-special of order $p^{2n+2}$, the genetic basis obtained in Proposition~\ref{genetic basis extra-special} consists of the group $S=P$, for which $N_P(S)/S$ is trivial, of $\frac{p^{2n+1}-1}{p-1}$ subgroups $S$ of index $p$ in $P$, for which $N_P(S)/S\cong C_p$, and the subgroup $S=Y$, for which $N_P(S)/S\cong Z(P)\cong C_{p^2}$. This gives
$$r=l_{\mathbb{R}}(P)=1+ \frac{p^{2n+1}-1}{p-1}\frac{p-1}{2}+\frac{p^2-p}{2}=\frac{p^{2n+1}+p^2-p+1}{2},$$ 
and
$$q=l_{\rac}(P)=1+\frac{p^{2n+1}-1}{p-1}+1=\frac{p^{2n+1}+2p-3}{p-1}.$$

This completes the proof.\endpf

To calculate $Cl_1(\Z P)$ we will need the following result.

\result{Lemma}
\label{segundo}
Let $p$ be an odd prime. 
Let $W$ be a vector space over $\F_p$ of finite dimension $k$, which is endowed with a bilinear, alternating form $b:W\times W\rightarrow \F_p$. Suppose that the rank of $b$ is not equal to 2. \par
For $x$ and $y$ in $W$, we still set $B_{x,\, y}=x^{p-1}y\in S^p(W)$. Then we have
\begin{displaymath}
S^p(W)=\langle B_{x,\,y}\mid x,\,y\in W \textrm{ s.t. }b(x,\, y)=0\rangle.
\end{displaymath}
\fresult
\pf
We will write $\mathcal{O}$ for $\langle B_{x,\,y}\mid x,\,y\in W \textrm{ s.t. } b(x,\, y)=0\rangle$. We will prove that $B_{x,y}\in \mathcal{O}$ for {\em any} $x,y\in W$, and by Lemma~\ref{Bxy generate}, it will follow that $S^p(W)=\mathcal{O}$.

Observe that, since $B_{x,y}\in\mathcal{O}$ for an $x,y\in W$ with $b(x,y)=0$, it follows from formula~\ref{the formula 2} that $x_1x_2\ldots x_{p-1}y\in\mathcal{O}$ for any elements $x_1,\ldots,x_{p-1}$ of $W$ such that $b(x_i,y)=0$ for $1\leq i\leq p-1$.

If $b=0$, there is nothing to prove. Otherwise let $x,y\in W$ such that $b(x,y)\neq 0$. Since $B_{x,\lambda y}=\lambda B_{x,y}$ for $\lambda\in\F_p$, to prove that $B_{x,y}\in\mathcal{O}$,  we can assume without loss of generality that $b(x,y)=1$, up to replacing $y$ by a suitable scalar multiple. If the restriction of $b$ to $\langle x,y\rangle^\perp$ was identically 0, then  $\langle x,y\rangle^\perp$ would be precisely the radical of $b$, so $b$ would have rank 2, contradicting our assumption. Hence we can find $z,t\in \langle x,y\rangle^\perp$ such that $b(z,t)\neq 0$, and up to replacing $t$ by some scalar multiple, we can assume that $b(z,t)=1$. \par
Now let $\alpha\in\F_p$, and set $u=\alpha x+t$ and $v=y+\alpha z$. Then 
$$b(u,v)=\alpha b(x,y)+\alpha^2b(x,z)+b(t,y)+\alpha b(t,z)=0,$$
since $b(x,y)=1=-b(t,z)$ and $b(x,z)=0=b(t,y)$.\par
It follows that $B_{u,v}\in\mathcal{O}$. But $B_{u,v}$ is equal to
\begin{equation}\label{mod O}
(\alpha x+t)^{p-1}(y+\alpha z)=\sum_{i=0}^{p-1}\binom{p-1}{i}\alpha^ix^iyt^{p-1-i} + \sum_{i=0}^{p-1}\binom{p-1}{i}\alpha^{i+1}t^{p-1-i}zx^i.
\end{equation}
By the observation at the beginning of the proof, the element $x^iyt^{p-1-i}$ is in $\mathcal{O}$ whenever $p-1-i>0$, since $b(x,t)=b(y,t)=b(t,t)=0$. Similarly $t^{p-1-i}zx^i\in\mathcal{O}$ if $i>0$, since $b(t,x)=b(z,x)=b(x,x)=0$. It follows that in (\ref{mod O}), the only elements possibly not in $\mathcal{O}$ correspond to $p-1-i=0$ in the first summation and to $i=0$ in the second. Hence
$$\alpha^{p-1}x^{p-1}y+\alpha t^{p-1}z\in\mathcal{O},$$
and this holds for any $\alpha\in\F_p$. For $\alpha=1$, this gives $x^{p-1}y+t^{p-1}z\in\mathcal{O}$, and for $\alpha=-1$, this gives $x^{p-1}y-t^{p-1}z\in\mathcal{O}$. It follows that $x^{p-1}y\in\mathcal{O}$, as was to be shown.\endpf

\begin{rem}{Remark} If the rank of $b$ is equal to 2, then the result of Lemma~\ref{segundo} is no longer true: for example in the non-degenerate case, that is when $W$ has dimension~2, saying that $b(x,y)=0$ for $x\neq 0$ is equivalent to saying that $y$ is a scalar multiple of~$x$. In this case $\mathcal{O}$ is the subspace of $S^p(W)$ generated by the elements $x^p$, for $x\in W$. So~$\mathcal{O}$ has dimension 2, whereas $S^p(W)$ has dimension $p+1$.
\end{rem}
\medskip
We now come to our main theorem, describing the structure of $Cl_1(\Z P)$ when $P$ is an extra-special or almost extra-special $p$-group for $p$ odd. We first recall that Oliver (\cite{bob} Example 7 page 16) showed that if $P$ is extra-special of order $p^3$, then  $Cl_1(\Z P)\cong(C_p)^{p-1}$, and that if $P$ is almost extra-special of order $p^4$, then $Cl_1(\Z P)\cong (C_p)^{(p^2+p-2)/2}$. Hence in what follows, we may assume that $P$ is an extra-special group of order at least $p^5$, or an almost extra-special $p$-group of order at least $p^6$.
\begin{rem}{Notation}\label{technique} Let $p$ be an odd prime and $n$ be a positive integer. Let $P$ be an extra-special $p$-group of order $p^{2n+1}$ or an almost extra-special $p$-group of order $p^{2n+2}$. Let $Z$ denote the center of $P$, and $N=P'\leq Z$ be the Frattini subgroup of $P$. {Let~$Y$ be a subgroup of $P$ of maximal order such that $Y\cap Z={\bf 1}$, as in Proposition~\ref{genetic basis extra-special}. Recall that $Y$ is elementary abelian. In any case, the group~$P$ can be written as a semidirect product $P=X\cdot Y$, where $X\geq Z$ is an abelian normal subgroup of $P$ with $C_P(X)=X$ and $Y\cap X={\bf 1}$:}
\begin{itemize}
\item If $P$ is extra-special of exponent $p$, the group $X$ is equal to $C\times X_0$, for some subgroup $X_0\cong (C_p)^n$ and $C=N=Z$.
\item  If $P$ is extra-special of exponent $p^2$, the group $X$ is equal to $C\times X_0$, for some subgroup $X_0\cong(C_p)^{n-1}$, some subgroup $C\cong C_{p^2}$, and $N=Z<C$. 
\item If $P$ is almost extra-special, then $X=C\times X_0$, for some subgroup $X_0\cong(C_p)^{n}$, and $N<Z=C\cong C_{p^2}$.
\end{itemize}
So in all cases we have $X=C\times X_0$, for some cyclic subgroup $C\geq Z\geq N$. \par
Moreover the subgroup~$Y$ is elementary abelian of order $p^n$.  It is maximal subject to the condition $Y\cap Z={\bf 1}$, so by Proposition~\ref{genetic basis extra-special}, we have a genetic basis of $P$ consisting of $P$ itself, its subgroups of index~$p$, and $Y$. The normalizer $N_P(Y)$ is equal to $Z\cdot Y$, so
$$\Gamma(P)\cong\Big(\prod_{[P:Q]=p}(P/Q)\Big)\times Z.$$


\end{rem}
\result{Theorem}\label{debut} Let $p$ be an odd prime, and let $P$ be an extra-special $p$-group of order at least $p^5$, or an almost extra-special $p$-group of order at least $p^6$. Let $N=P'$ be the Frattini subgroup of $P$, and $Z$ be the center of $P$. Then there is a split sequence of abelian groups
$$\xymatrix{
0\ar[r]&K\ar[r]&Cl_1(\Z P)\ar[rr]^-{\dfl_{P/N}^P}&&Cl_1\big(\Z (P/N)\big)\ar[r]&0,\\
}
$$
where $K$ is cyclic, 
isomorphic to a quotient of $Z$. In particular $Cl_1(\Z P)$ is isomorphic to $K\times (C_p)^M$, where $K$ is cyclic of order dividing~$p^2${\rouge, and $M=\frac{p^{k-1}-1}{p-1}-\binom{p+k-2}{p}$ if $|P|=p^k$}.
\fresult
\pf The product $\prod_{[P:Q]=p}\limits(P/Q)$ identifies with $\Gamma(P/N)$, and we have a surjective projection map $\dfl_{P/N}^P:\Gamma(P)\to \Gamma(P/N)$, with kernel isomorphic to $Z$. By Proposition~\ref{ladef}, this map induces a surjective deflation map 
$$\dfl_{P/N}^P:Cl_1(\Z P)=\Gamma(P)/\mathcal{R}\to Cl_1\big(\Z (P/N)\big)=\Gamma(P/N)/\sur{\mathcal{R}},$$
where $\mathcal{R}$ is the subgroup of $\Gamma(P)$ generated by the elements $u_{H,g}$ introduced in Theorem~\ref{elmejor}, and $H$ is a cyclic subgroup of $P$ with $g\in C_P(H)$. Similarly $\sur{\mathcal{R}}$ is the corresponding subgroup of $\Gamma(P/N)$ generated by the elements $u_{F,c}$, where $F$ is a cyclic subgroup of $P/N$ and $c\in P/N$ (we always have $c\in C_{P/N}(F)$, as $P/N$ is abelian). \par
The proof of Theorem~\ref{ladef} shows that $\dfl_{P/N}^P(u_{H,g})=u_{HN/N,gN}$. Conversely, if $F$ is a cyclic subgroup of $P/N$, generated by $f$,  and if $c\in P/N$, then there exists a pair $(H,g)$ of a cyclic subgroup $H$ of $P$ and an element $g\in C_P(H)$ such that $HN/N=F$ and $gN=c$ if and only if $b(f,c)=0$, where $b$ is the bilinear alternating form on $P/N$ with values in $\F_p$ induced by taking commutators in $P$. Our assumptions on $P$ imply that the rank of $b$ is not 2, so we can apply Lemma~\ref{segundo}. This shows that the subspace of $\Gamma(P/N)$ generated by the elements  $u_{F,c}$, where $F=\langle{f}\rangle$ for $f\in P/N$, and $c\in P/N$ such that $b(f,c)=0$, generate $S^p(P/N)=\sur{\mathcal{R}}$, by Lemma~\ref{primero}. It follows that $\dfl_{P/N}^P$ induces a {\em surjective} map $e:\mathcal{R}\to\sur{\mathcal{R}}$. Let $L$ denote the kernel of this map. We have a commutative diagram with exact rows
$$\xymatrix{
0\ar[r]&L\ar[d]^l\ar[r]&\mathcal{R}\ar[d]^-i\ar[rr]^-{e}&&\sur{\mathcal{R}}\ar[d]^-j\ar[r]&0\\
0\ar[r]&Z\ar[r]&\Gamma(P)\ar[rr]^-{\dfl_{P/N}^P}&&\Gamma(P/N)\ar[r]&0\\
}
$$
where the vertical maps $i$ and $j$ are the inclusion maps. The Snake's Lemma now shows that the map $l$ is injective, and moreover we have an exact sequence of cokernels
\begin{equation}\label{split}
\xymatrix{
0\ar[r]&K\ar[r]&Cl_1(\Z P)\ar[rr]^-{\dfl_{P/N}^P}&&Cl_1\big(\Z (P/N)\big)\ar[r]&0,
}
\end{equation}
\!where $K=Z/l(L)$. Since the kernel $Z$ of $\dfl_{P/N}^P$ corresponds to the component of $\Gamma(P)$ indexed by $Y$, the image $l(L)$ is generated by the components $w_Y\in N_P(Y)/Y\cong Z$ of vectors $w$ in the kernel of the deflation map $e:\mathcal{R}\to\sur{\mathcal{R}}$. Hence $l(L)$ is a subgroup of the group generated by all the elements $u_{H,g,Y}$, where $H$ is a cyclic subgroup of $P$ and $g\in C_P(H)$.\par
It remains to see that the exact sequence~\ref{split} splits. To see this, consider the completed diagram
$$\xymatrix{
&0\ar[d]&0\ar[d]&&0\ar[d]\\
0\ar[r]&L\ar[d]^l\ar[r]&\mathcal{R}\ar[d]^-i\ar[rr]^-e&&\sur{\mathcal{R}}\ar[d]^-j\ar[r]&0\\
0\ar[r]&Z\ar[d]\ar[r]&\Gamma(P)\ar[d]^-a\ar[rr]^-{c}&&\Gamma(P/N)\ar@/^2ex/@{-->}[ll]^-t\ar[d]^-b\ar[r]&0\\
0\ar[r]&K\ar[d]\ar[r]&Cl_1(\Z P)\ar[d]\ar[rr]^-{d}&&Cl_1\big(\Z (P/N)\big)\ar@/^1ex/@{-->}[u]^-s\ar[d]\ar[r]&0,\\
&0&0&&0\\
}
$$
where $a$ and $b$ are the projection maps, and $c$, $d$ are the respective deflation maps, and $e$ is the restriction of $c$ to $\mathcal{R}$. The map $b$ is split surjective, because $\Gamma(P/N)$ and $Cl_1\big(\Z(P/N)\big)$ are both elementary abelian.  Similarly, the map $c$ is split surjective by construction. Let $s$ be a section of $b$ and let $t:\Gamma(P/N)\to \Gamma(P)$ be a section of $c$. Then 
$$d\circ a\circ t\circ s=b\circ c\circ t\circ s=b\circ s={\rm Id},$$
so the map $a\circ t\circ s$ is a section of $d$. To complete the proof, observe that $K$ is a isomorphic to a quotient of $Z$, and that $Z$ is cyclic of order $p$ or $p^2$.\endpf
\vspace{-2ex}
\begin{rem}{Remark}\label{generators} It follows from this proof that to obtain a minimal set $\mathcal{M}$ of generators of $Cl_1(\Z P)$, we can take a subset $(v_S)_{S\in\mathcal{S}}$ of the canonical basis of $\Gamma(P/N)$ (as $\F_p$-vector space), corresponding to a set $\mathcal{S}$ of subgroups of index $p$ of $P$, which maps by $b$ to a basis of $Cl_1\big(\Z (P/N)\big)$. Then $v_S$ is a generator {of the component $P/S\cong (P/N)\big/(S/N)$ of $\Gamma(P)$, for $S\in\mathcal{S}$. If $K$ is trivial\footnote{This will occur only in {\rouge Assertion 1 of Proposition~\ref{fait chier}} below, i.e. when $P$ is extra-special of order $p^5$ and exponent $p^2$}, we set $\mathcal{T}=\mathcal{S}$. Otherwise, we also choose a generator $v_Y$ of $Z\cong N_P(Y)/Y$, and we set $\mathcal{T}=\mathcal{S}\sqcup\{Y\}$. Then we take for $\mathcal{M}$ the image in $Cl_1(\Z P)$ by the projection map $a$ of the set $\{v_S\mid S\in\mathcal{T}\}$.}\par
As we can see from this procedure, the splitting of the exact sequence of Theorem~\ref{debut} is {\em not} canonical: first we have made a choice for the genetic subgroup $Y$, and then we also make a choice of a set $\mathcal{S}$ and of generators~$v_S$, for $S\in\mathcal{T}$.\nopagebreak
\end{rem}\vspace{1ex}

{We now come to the proof of Theorem A, which we split in two parts: first the ``generic" case, when the group $P$ is large enough, and then two special cases of ``small" groups of order $p^5$ and~$p^6$. The generic case is} {the following one.}
{
\result{Theorem}\label{principal} Let $p$ be an odd prime, and let $P$ be a $p$-group of order $p^k$ which is either: \begin{enumerate}
\item extra-special of exponent $p$ with $k\geq 5$, or 
\item extra-special of exponent $p^2$ with $k\geq 7$, or 
\item almost extra-special  with $k\geq 8$. 
\end{enumerate}
Then {$Cl_1(\Z P)$ is (non-canonically)} isomorphic to $Z\times (C_p)^M$, where $Z$ is the center of~$P$ (so $Z\cong C_p$ in cases 1 and 2, and $Z\cong C_{p^2}$ in case 3), and $M=\frac{p^{k-1}-1}{p-1}-\binom{p+k-2}{p}$.
\fresult
\pf {We keep Notation~\ref{technique} throughout.} We know from Theorem~\ref{debut} that $Cl_1(\Z P)\cong K\times (C_p)^M$, where $M=\frac{p^{k-1}-1}{p-1}-\binom{p+k-2}{p}$. \par
Moreover the group $K$ is isomorphic to the quotient of $Z$ by the subgroup generated by the elements $w_Y$, where $w$ is an element in the kernel of the deflation map $e:\mathcal{R}\to \sur{\mathcal{R}}$. Such an element $w_Y$ is a product of elements $u_{H,g,Y}$, for some pairs $(H,g)$ of a cyclic subgroup $H$ of $P$, and an element $g\in C_P(H)$. We start by computing $u_{H,g,Y}$ for such a pair $(H,g)$, using formula~\ref{formule2}. We will show that $u_{H,g,Y}=1$ in {the cases {of} the statement of the theorem}, 
and it will follow that $K\cong Z$. We have
\begin{displaymath}
u_{H,g,Y}=\prod_{x\in D}\overline{l_{g,\, x}}
\end{displaymath}
where $D$ is a chosen set of representatives of those double cosets $H\langle g\rangle xN_P(Y)$ in $P$ for which $H^x\cap N_P(Y)\leq Y$. For each $x\in D$, if $m_x$ is the index of $\langle g\rangle\cap {\big(}H\cdot{^xN_P(Y)}{\big)}$ in $\langle g\rangle$, the element $g^{m_x}$ can be written $g^{m_x}=h_x\cdot {^xl_{g,x}}$, for $h_x\in H$ and $l_{g,x}\in N_P(Y)$. Since $N_P(Y)=ZY$ is a normal subgroup of $p$, by Lemma~\ref{deboma}, the group $H\cdot{^xN_P(Y)}$ and the integer $m_x$ do not depend on $x$. Let $m$ denote this integer. Since $g^p\in Z$ for any $g\in P$, we have $m=1$ if $g\in HZY$, and $m=p$ otherwise. \par
Since $H\langle g\rangle ZY$ is also a normal subgroup of $P$, we have $H\langle g\rangle tZY=H\langle g\rangle ZYt=tH\langle g\rangle ZY$, for any $t\in P$. Hence $D$ is a subset of a set of representatives of $P/H\langle g\rangle ZY$. If $D$ is empty, we have $u_{H,g,Y}=1$, so we can assume $D\neq\emptyset$. If $H^t\cap ZY\leq Y$, then $H^t\cap Z\leq Y\cap Z={\bf 1}$. Then $H\cap Z={\bf 1}$, so $H$ has order 1 or $p$, since $h^p\in Z$ for any $h\in P$. 
\begin{itemize}
\item If $H\cap ZY={\bf 1}$ then $H^t\cap ZY=(H\cap ZY)^t={\bf 1}$ for any $t\in P$, so our set $D$ is a set of representatives of $P/H\langle g\rangle ZY$. 
\item If $H\cap ZY\neq {\bf 1}$, then $H\leq ZY$. If there exists $t_0\in P$ with $H^{t_0}\cap ZY\leq Y$, we have $H^{t_0}\leq Y$, and up to replacing $H$ by $H^{t_0}$, we can assume $H\leq Y$. Then for $t\in P$, we have $H^t\cap ZY=H^t$, and $H^t\leq Y$ implies $[H,H^t]\leq Z\cap Y={\bf 1}$. In other words $H^t\cap ZY\leq Y$ if and only if $t\in C_P(H)$. Here $D$ is a set of representatives of $C_P(H)/H\langle g\rangle ZY=C_P(H)/\langle g\rangle ZY$.
\end{itemize}
So in each case, there are subgroups $M$ and $Q$ of $P$ with $ZY\leq M\leq Q$ such that $D$ is a set of representatives of $Q/M$.\par

We have $g^m=h{\cdot}z{\cdot}y$, for some $h\in H$, $z\in Z$ and $y\in Y${. O}bserve that $z=g^p$ and $h=y=1$ if $m=p$. Then for $x\in D$, we have $g^m=h{\cdot}z{\cdot}^x(y^x)=h{\cdot}^x(z[y,x]y)$, so with the identification $N_P(Y)/Y\cong Z$, we have $l_{g,x}=z[y,x]$. It follows that $u_{H,g,Y}=\prod_{x\in D}\limits\big(z[y,x]\big)=z^{|D|}[y,\prod_{x\in D}\limits x]$. Now the commutator $[y,\prod_{x\in D}\limits x]$ only depends of the images of $y$ and $\prod_{x\in D}\limits x$ in the elementary abelian group $P/N$. As $M\geq N$, the map $x\mapsto \sur{x}=xN$ is a bijection from $Q/M$ to $\sur{Q}/\sur{M}$, where $\sur{Q}=Q/N$ and $\sur{M}=M/N$. The subgroup $\sur{M}$ admits a supplement $W$ in the elementary abelian group $\sur{Q}$, and $W$ is a set of representatives of $\sur{Q}/\sur{M}$. Since $\sum_{w\in W}\limits w=0$ as $p$ is odd, it follows that $\prod_{x\in D}\limits x$ maps to 0 in $P/N$, that is $\prod_{x\in D}\limits x\in N\leq Z$, and $[y,\prod_{x\in D}\limits x]=1$. This gives finally
\begin{equation}\label{uHgY}
u_{H,g,Y}=\left\{\begin{array}{cl}z^{|D|}&\hbox{if}\;g\in HzY\;\hbox{for some}\;z\in Z,\\g^{p|D|}&\hbox{otherwise.}\end{array}\right.
\end{equation}
\mpn
{$\bullet$} If $P$ is extra-special of order $p^{2l+1}$, with $l\geq 2$, then $Y$ has order $p^l$ and $Z$ has order~$p$.
\begin{itemize}
\item[-] {As we said,} if $g\in HzY$, for $z\in Z$, {we have} $u_{H,g,Y}=z^{|D|}$. If $H\nleq ZY$, then 
$$|D|=|P:H\langle g\rangle ZY|=|P:HZY|\geq p^{2l+1}/(p{\cdot}p{\cdot}p^l)=p^{l-1}.$$
Now if we assume $H\leq Y$, then $g\in HZY=ZY$, and
$$|D|=|C_P(H):\langle g\rangle ZY|=[C_P(H):ZY|\geq p^{2l}/p^{l+1}=p^{l-1}$$
again. For $l\geq 2$, this is a multiple of $p=|Z|$, and $u_{H,g,Y}=1$.
\item[-] If $g\notin HZY$, then $u_{H,g,Y}=g^{p|D|}$. So $u_{H,g,Y}=1$ unless $g$ has order $p^2$, and then $g^p$ generates $Z$. If $H\nleq ZY$, then 
$$|D|=|P:H\langle g\rangle ZY|=|P:H\langle g\rangle Y|\geq p^{2l+1}/p{\cdot}p^2{\cdot} p^l=p^{l-2}.$$
{On the other hand,} if $H\leq Y$, then 
$$|D|=|C_P(H):\langle g\rangle ZY|=|C_P(H):\langle g\rangle Y|\geq p^{2l}/p^2{\cdot}p^l=p^{l-2}$$
again.  Hence $|D|$ is a multiple of $p$ if $l\geq3$, so $u_{H,g,Y}=1$.
\end{itemize}
So $u_{H,g,Y}=1$ if $P$ is extra-special, unless possibly if $P$ has exponent $p^2$ and order~$p^5$.
\mpn{$\bullet$} If $P$ is almost extra-special of order $p^{2l+2}$, with $l\geq 2$, then $Y$ has order $p^l$ and $Z$ is cyclic of order $p^2$.\par
\begin{itemize}
\item[-] If $g\in HzY$, for $z\in Z$, then $u_{H,g,Y}=z^{|D|}$. If $H\nleq ZY$, we have 
$$|D|=|P:H\langle g\rangle ZY|=|P:HZY|\geq p^{2l+2}/(p{\cdot}p^2{\cdot}p^l)=p^{l-1}.$$
If we now assume $H\leq Y$, then  $g\in HZY=ZY$, and
$$|D|=|C_P(H):\langle g\rangle ZY|=|C_P(H):ZY|\geq p^{2l+1}/p^2{\cdot} p^l=p^{l-1}$$
again. Hence then $|D|$ is a multiple of $p^2$ if $l\geq 3$, so $u_{H,g,Y}=1$. 
\item[-] If $g\notin HZY$, then $u_{H,g,Y}=g^{p|D|}$. If $H\nleq ZY$, since $\langle g\rangle Z$ has order at most $p^3$ as $g^p\in Z$, we have 
$$|D|=|P:H\langle g\rangle ZY|\geq p^{2l+2}/p{\cdot}p^3{\cdot}p^l=p^{l-2},$$
If we now assume $H\leq Y$, then 
$$|D|=|C_P(H):\langle g\rangle ZY|\geq p^{2l+1}/p^3{\cdot}p^l=p^{l-2}$$
again. Hence then $|D|$ is a multiple of $p$ if $l\geq 3$, and $u_{H,g,Y}=g^{p|D|}=1$.
\end{itemize}
So $u_{H,g,Y}=1$ if $P$ is almost extra-special, unless possibly if $P$ has order~$p^6$.\endpf

So we are left with the special cases where $P$ is either extra-special of order $p^5$ and exponent $p^2$, or almost extra-special of order $p^6$. In both cases, we will use the {next} lemma, where the notation is again as in \ref{technique}{.}
\result{Lemma}\label{abU} {\rouge Let $P$ be extra-special of order $p^5$ and exponent $p^2$, or almost extra-special of order $p^6$. Then t}here exist $a\in X_0$ and $b\in Y$ with the following properties:
\begin{enumerate}
\item The elements $a$ and $b$ both have order $p$, and centralize $C$.
\item The group $U=\langle a,b\rangle$ is extra-special of order $p^3$ and exponent $p$. Its center is equal to $N$.
\item For $\alpha\in\{0,\ldots,p-1\}$, set $H_\alpha=\langle ab^\alpha\rangle$. Also set $H_\infty=\langle b\rangle$, and denote by $\mathbb{L}$ the set $\{0,\ldots,p-1,\infty\}$. Then $\langle H_\alpha,H_\beta\rangle=U$ for any pair $(\alpha,\beta)$ of distinct elements of $\mathbb{L}$, and the map $\alpha\in\mathbb{L}\mapsto H_\alpha N$ is a bijection from $\mathbb{L}$ to the set of subgroups of order $p^2$ of $U$.
\end{enumerate}
\fresult
\pf In both cases $Y\cong (C_p)^2$. Since $[C,Y]\leq Z\leq C$, the group $Y$ normalizes~$C$, and $C_Y(C)\neq {\bf 1}$ since the automorphism group of $C$ is cyclic. We choose a non trivial element $b$ of $C_Y(C)$. Then $b$ does not centralize $X_0$, for otherwise $b$ centralizes $CX_0=X$, contradicting $C_P(X)=X$. Hence we can choose $a\in X_0$ such that $[a,b]\neq 1$. Then the group $U=\langle a,b\rangle$ is extra-special of order $p^3$ and exponent $p$, as $a$ and $b$ have order $p$ and $[a,b]$ generates $N$. In particular $N$ is equal to the center of $U$.\par
The groups $H_\alpha$, for $\alpha\in\mathbb{L}$, all have order $p$, and are different from $N$. So the groups $H_\alpha N$ are subgroups of order $p^2$ of $U$. Clearly $\langle H_\alpha,H_\beta\rangle=U$ for any distinct element $\alpha, \beta$ of $\mathbb{L}$, hence the subgroups $H_\alpha N$, for $\alpha\in\mathbb{L}$, are all distinct. This completes the proof, as $U$ has exactly $p+1=|\mathbb{L}|$ subgroups of order $p^2$.\endpf \par

{We now proceed with the proof of the above mentioned two special cases.}
{\rouge
\result{Proposition} \label{fait chier}\begin{enumerate}
\item Let $P$ be extra-special of order $p^5$ and exponent $p^2$. Then the group $K$ of Theorem \ref{debut} is trivial.
\item Let $P$ be almost extra-special of order at least $p^6$. Then the group $K$ of Theorem \ref{debut} has order $p$.
\end{enumerate}
\fresult
\pf We keep the notation of Lemma~\ref{abU}.\par\noindent
1) Suppose first that $P$ is extra-special of order $p^5$ and exponent $p^2$.} Then we have $X=C\times X_0$, where $C\cong C_{p^2}$ contains the center $Z=N$ of order $p$ of $P$, and $X_0\cong C_p$. On the other hand $Y\cong(C_p)^2$ in this case. We will build an element $w$ in the kernel of $e:\mathcal{R}\to\sur{\mathcal{R}}$ with non trivial component $w_Y\in Z$. Then $w_Y$ will generate~$Z$, so $K$ will be trivial in this case, as stated. \par
In order to build $w$, we first fix a generator $g$ of $C$. Next we consider {\rouge the set $\mathbb{L}$ and the subgroups $H_\alpha$ of $C_P(C)=C_P(g)$, for $\alpha\in\mathbb{L}$, introduced in Lemma~\ref{abU}. We now set}
$$w=(u_{{\bf 1},g})^{-1}\prod_{\alpha\in \mathbb{L}}u_{H_\alpha,g}.$$
For an arbitrary cyclic subgroup $H$ of $C_P(g)$, the component $u_{H,g,S}$ of $u_{H,g}$ at a genetic subgroup $S$ of $P$ is given by Formula~\ref{formule2}. When $|P:S|=p$, this formula comes down to
$$u_{H,g,S}=\left\{\begin{array}{cl}1&\hbox{if}\;H\nleq S,\\\sur{g}&\hbox{otherwise},\end{array}\right.$$
where $\sur{g}$ is the image of $g$ in $P/S$. It follows that the component $w_S$ of $w$ at $S$ is equal to $(\sur{g})^{-1+n_S}$, where $n_S=|\{\alpha\in\mathbb{L}\mid H_\alpha\leq S\}|$. Then either $S\geq U=\langle a,b\rangle$, and $S\geq H_\alpha$ for all $\alpha\in\mathbb{L}$, so $n_S=p+1$, and $w_S=(\sur{g})^p=1$, since $P/S$ has order $p$; or $S\ngeq U$, and then $S\cap U$ is a subgroup of index $p$ of $U$. So $|S\cap U|=p^2$, and $S\cap U=H_\alpha Z$ for a unique $\alpha\in\mathbb{L}$. Hence $H_\alpha\leq S$, and $H_\beta\nleq S$ for $\beta\in \mathbb{L}-\{\alpha\}$. It follows that $n_S=1$ in this case, so $w_S=(\sur{g})^{-1+1}=1$ again. This shows that $w$ is in the kernel of the map $e:\mathcal{R}\to\sur{\mathcal{R}}$.\par
We now compute the component $w_Y$ of $w$ indexed by $Y$, using~\ref{uHgY}: we consider a cyclic subgroup $H$ of  order at most $p$ of $C_P(g)$. Then $g\notin HZY$, because $HZY$ has exponent $p$. Then $u_{H,g,Y}=g^{p|D|}$, where $D$ is a set of representatives of those double cosets $H\langle g\rangle tZY$ such that $H^t\cap ZY\leq Y$. \par
We know moreover that $|D|=|P:H\langle g\rangle ZY|=|P:HCY|$ if $H\nleq ZY$, since $Z\leq C=\langle g\rangle$, and that $|D|=|C_P(H):\langle g\rangle ZY|=|C_P(H):CY|$ if $H\leq Y$.\par
If $H={\bf 1}$, then $|D|=|P|/(|C||Y|)=p^5/(p^2{\cdot}p^2)=p$, so $u_{H,g,Y}=g^{p^2}=1$. If $H=H_\alpha$, for $\alpha\in\mathbb{L}-\{\infty\}$, then $|D|=|P:H_\alpha CY|=p^5/(p{\cdot}p^2{\cdot}p^2)=1$, {for} $H_\alpha\nleq CY$. Then $u_{H_\alpha,g,Y}=g^p$. Finally, if $H=H_\infty$, we have $|D|=|C_P(H_\infty):CY|=p^4/(p^2{\cdot}p^2)=1$, so $u_{H_\alpha,g,Y}=g^p$ again.\par

This gives $w_Y=1^{-1}\prod_{\alpha\in \mathbb{L}}\limits g^p=g^{p(p+1)}=g^p$, and $g^p$ is a generator of $Z$. So we have indeed built an element $w$ in the kernel of $e:\mathcal{R}\to \sur{\mathcal{R}}$ such that the component $w_Y$ generates $Z$. It follows that the group $K$ is trivial in this case, as claimed.\medskip\par\noindent
{\rouge 2)} Suppose now that $P$ is almost extra-special of order $p^6$. Then we have $C=Z\cong C_{p^2}>P'=N$. In this case, let $g\in P$ and $H$ be a subgroup of order 1 or $p$ of $C_P(H)$, Equation~\ref{uHgY} shows that $u_{H,g,Y}$ is equal to $(v_{H,g})^{|D|}$ for some $v_{H,g}\in Z$, where $D$ is a set of representatives of double cosets $H\langle g\rangle tZY$ such that $H^t\cap ZY\leq Y$. We claim that $u_{H,g,Y}\in N${. I}ndeed, if $g\notin HZY$, then $v_{H,g}=g^p\in N${. So we can assume that} $g\in HZY$. Then if $H\nleq ZY$, we have $|D|=|P:H\langle g\rangle ZY|=|P:HZY|\geq p^6/(p{\cdot}p^2{\cdot}p^2)=p$. If $H\leq Y$, then $g\in HZY=ZY$, and $|D|=|C_P(H):ZY|\geq p^5/(p^2{\cdot}p^2)=p$. In both cases $|D|$ is a multiple of $p$, so $u_{H,g,Y}=(v_{H,g})^{|D|}\in N$ again, as claimed. \par
As in {\rouge the proof of Assertion 1}, we will build an element $w$ in the kernel of $e:\mathcal{R}\to\sur{\mathcal{R}}$ such that $w_Y$ generates $N$. {This} will show that $K\cong Z/N\cong C_p$, as claimed. We first choose a generator $g$ of $Z=C$. The formal definition of $w$ is then the same as in {\rouge the proof of Assertion 1}:
$$w=(u_{{\bf 1},g})^{-1}\prod_{\alpha\in \mathbb{L}}u_{H_\alpha,g}.$$
The computation of the component $w_S$ for a subgroup $S$ of index $p$ of $P$ is also very similar to what we did for Assertion 1: for a subgroup $S$ of index $p$ of $P$, we have $w_{S}=(\sur{g})^{-1+n_S}$, where $n_{S}=|\{\alpha\in\mathbb{L}\mid H_\alpha\leq S\}|$. We have $n_S=1$ or $n_S=p+1$, thus $w_S=1$, and $w$ lies in the kernel of $e:\mathcal{R}\to\sur{\mathcal{R}}$.\par

As for the component $w_Y$ of $w$, we observe that $u_{H,g,Y}=g^{|D|}$ since $g\in Z$, where $D$ is as above. If $H={\bf 1}$, then $|D|=|P:ZY|=p^2$, and $u_{{\bf 1},g,Y}=g^{|P:ZY|}=g^{p^2}=1$. For $\alpha\in\mathbb{L}-\{\infty\}$, we have $|D|=|P:H_\alpha \langle g\rangle ZY|=|P/H_\alpha ZY|=p^6/p{\cdot}p^2{\cdot}p^2=p$, so $u_{H_\alpha,g,Y}=g^p$. Finally, for $H=H_\infty\leq Y$, we have $|D|=|C_P(H_\infty):\langle g\rangle ZY|=|C_P(H_\infty):ZY|=p^5/(p^2{\cdot}p^2)=p$, so $u_{H_\infty,g,Y}=g^p$ again.\par

It follows that $w_Y=1^{-1}(g^p)^{p+1}=g^{p^2+p}=g^p$, so $w_Y$ generates $N$, as was to be shown.\endpf
}

{
\centerline{\rule{5ex}{.1ex}}
\begin{flushleft}
Serge Bouc, CNRS-LAMFA, 33 rue St Leu, 80039, Amiens, France.\\
{\tt serge.bouc@u-picardie.fr}\vspace{1ex}\\
Nadia Romero, DEMAT, UGTO, Jalisco s/n, Mineral de Valenciana, 36240, Guanajuato, Gto., Mexico.\\
{\tt nadia.romero@ugto.mx}
\end{flushleft}
}
\end{document}